\documentclass{article}
\usepackage{amsfonts}
\usepackage{amsmath}
\usepackage{amssymb}
\usepackage{cite}
\usepackage{fullpage}
\usepackage{graphicx}
\usepackage{dsfont}
\def\bR{{\mathbb{R}}}
\def\bE{{\mathbb{E}}}
\def\dx{{\Delta x}}
\def\dt{{\Delta t}}
\def\eps{\varepsilon}

\def\calF{{\cal F}}

\def\one{\mathds{1}}

\def\labs{\left\lvert}
\def\rabs{\right\rvert}

\title{Non-parametric estimation of Stochastic Differential Equations from stationary time-series}

\date{\today}

\author{Xi Chen\thanks{University of Houston, Department of Mathematics, 
victoryx@math.uh.edu} \and
Ilya Timofeyev\thanks{University of Houston, Department of Mathematics,      
ilya@math.uh.edu}
}

\begin{document}
\maketitle
\begin{abstract}
We study efficiency of non-parametric estimation of diffusions (stochastic differential equations driven by Brownian motion) from long stationary trajectories.
First, we introduce estimators based on conditional expectation which is motivated by the definition of drift and diffusion coefficients. 
These estimators involve time- and space-discretization parameters for computing expected values from discretely-sampled stationary data.
Next, we analyze consistency and mean squared error of these estimators depending on 
computational parameters.
We derive relationships between the number of observational points, time- and space-discretization parameters
in order to achieve the optimal speed of convergence and minimize computational complexity. 
We illustrate our approach with numerical simulations. 

\end{abstract}

\section{Introduction}
\label{sec:intro}

Recently, there has been a significant increase in the amount of available observational data. 
Various areas, such as 
biology, geosciences, social science, etc. provide large datasets which need to be analyzed. 
In particular, it is often necessary to fit an empirical model using available stationary data with the goal of forecasting 
future values or generating trajectories with similar statistical properties. Such examples for instance often 
arise in turbulence (e.g. \cite{sura03,suba02,kkg05,nn2,dta12,dat12,mtv4,mtv2,zdat18,resseguier20}), 
reduced modeling of nonlinear dynamics (e.g. \cite{crove06a,klin,nn1,mch,nn3,mtv6}), 
and biology (e.g. \cite{clementi1,clementi2,schutte1}).
This is a very active area of research with many publications including results on parametric and non-parametric estimation of autoregressive processes and stochastic differential equations.

In this paper we elucidate how to optimally select computational parameters for non-parametric estimation
of the drift and diffusion coefficients in stochastic differential equations from discretely sampled stationary data.
Compared to parametric techniques, non-parametric approaches are, typically, more computationally challenging, 
but exhibit more flexibility, since non-parametric estimation does not rely on a particular functional form of the drift and diffusion coefficients. 
Therefore, we provide guidelines for reducing the computational complexity if the non-parametric estimation based on conditional expectations while maintaining the 
accuracy of drift and diffusion estimators.

Recently,  several authors explored non-parametric estimation approach based on
conditional expectations \cite{suba02,clementi1,lele10,schutte1}. 
This estimation technique relies on discrete analogs of conditional expectations
which are used to define the drift and diffusion coefficients in stochastic differential 
equations driven by Brownian motion \cite{ga85,book:oksendal}.
Since the estimation is non-parametric, it does not require any a-priori anzatz about 
the functional form of the drift and diffusion coefficients. 
Therefore, this estimation technique is quite general and can be applied to numerical and
experimental data without restricting the drift and diffusion coefficients to a particular form (e.g. additive noise only).
On the other hand, similar to other non-parametric techniques, 
the conditional expectation estimation requires a substantial amount of data. 
Therefore, it is essential to address the computational efficiency of this estimation technique
in practical situations.  In this paper we analyze the relationship between the space- and time-discretization parameters and derive  an explicit criteria for selecting these computational parameters in order to make this approach computationally efficient. 

The rest of the paper is organized as follows. In sections \ref{sec:sec1} and 
\ref{sec:sec3} we introduce background material for stochastic differential equations and
It\^o-Taylor Expansions, respectively. In sections \ref{sec:sec4} and \ref{sec:sec5}
we analyze the bias and the mean-squared error of the drift and diffusion estimators.
Section \ref{sec:sec5} contains main analytical results with details of analytical calculations presented in the appendix. Numerical results are presented in section \ref{sec:sec6}. Finally, we summarize our results in section \ref{sec:conc}.

\section{Non-parametric estimation of drift and diffusion}
\label{sec:sec1}

In this paper we consider one-dimensional stochastic differential equations (SDEs) 
driven by Brownian motion
\begin{equation}
\label{sde}
dX_t = A(X_t) dt + D(X_t) dW_t
\end{equation}
where, for simplicity, $X_t \in \bR^1$ and $W_t$ is 1-dimensional Brownian motion.
Our analysis can be generalized for $X_t \in \bR^n$ and $W_t \in \bR^m$, but mathematical expressions
become cumbersome and more difficult to read. Analysis of the one-dimensional case provides sufficient
guidelines for understanding behavior of estimators and optimal selection of computational parameters.
We assume that neither drift $A(X_t)$ nor diffusion coefficient $D(X_t)$ depend 
explicitly on time and, moreover, the SDE in \eqref{sde} has a unique stationary distribution
$\rho(x)$ such that $L^{FP} \rho(x) = 0$ where $L^{FP}$ is the Fokker-Planck operator given by
\begin{equation}
\label{Lfp}
L^{FP} = -\frac{\partial}{\partial x}A(x) + \frac{1}{2}\frac{\partial^2}{\partial x^2} D^2(x).
\end{equation}
Then it is known from the theory of parabolic equations \cite{ga85,book:risken} that under appropriate conditions distribution of $X_t$ converges to $\rho(x)$ as $t\to\infty$.
We also assume that $A(x)$ and $D(x)$ are sufficiently differentiable with finite derivatives 
since we're using It\^o-Taylor expansions in this paper.
This implies that $A(x)$ and $D^2(x)$ are uniformly Lipschitz on bounded intervals.
In addition, drift and diffusion coefficients of the SDE in \eqref{sde} can be 
defined as conditional expectations 
\cite{ga85,book:oksendal}
\begin{eqnarray}
A(x) &=&  \lim_{\dt \to 0}\frac{1}{\dt} \mathbb{E} \left[ X_{\dt} - x \big| X_0 = x \right],
\label{Ax} \\
D^2(x) &=& \lim_{\dt \to 0}\frac{1}{\dt} \mathbb{E} \left[(X_{\dt} - x)^2 \big| X_0  = x\right].
\label{Bx}
\end{eqnarray}
%

\subsection{Estimators for drift and diffusion from stationary time-series}
\label{sec2.1}

In several papers (e.g. \cite{suba02,clementi1,lele10,schutte1}) 
authors used definitions \eqref{Ax} and \eqref{Bx} to develop numerical approaches for 
estimating the drift and diffusion coefficients from stationary trajectories. 
To develop non-parametric estimators for the drift and diffusion coefficients based on conditional expectation 
in \eqref{Ax}, \eqref{Bx} we consider the following setup.

Assume that the available data are sampled from a stationary time-series of $X_t$ with a uniform 
time-step, $\dt$,  i.e. the available data are $\{U_k = X_{k\dt}, \, k=1,\ldots,N\}$. 
To develop practical approach for non-parametric estimation of drift and diffusion coefficients from 
such data we need to introduce estimators conditioned on an interval and not on a particular value since
it is extremely unlikely that for any given $x$ we can find any $k$ such that $U_k = x$,  i.e. the time-series are
unlikely to contain values exactly equal to $x$ for any given $x$. Even if we try to estimate the drift and diffusion at
$x=U_1$, the probability that $U_k = x$ for $k>1$ is zero. Moreover, in practical situations the goal is to estimate the
drift and diffusion at many values of $x$ (possibly discrete with a certain space-step). 
Therefore, we introduce a discrete uniform mesh in state-space, $x_k$ for $k=1,\ldots,K$ with
$x_{k+1} - x_k = \dx.$ Points $x_k$ represent centers of bins $Bin_k = [x_k - \dx/2, x_k + \dx/2]$
for computing analogs of the expected values in \eqref{Ax}, \eqref{Bx} numerically.
In practice space-discretization does not have to be uniform, but varying $\dx$ 
does not affect our results since our error analysis is performed for each bin separately.
Thus, we introduce discrete estimators 
for $A(x_k)$ and $D^2(x_k)$ as follows
\begin{eqnarray}
\hat{A}(x_k) &=&\frac{1}{M\dt}  \sum\limits_{j \in M_k} (X_{t_j + \dt} - X_{t_j}),
\label{Ahat} \\
\hat{D}^2(x_k) &=&\frac{1}{M\dt}  \sum\limits_{j \in M_k} \left(X_{t_j + \dt} - X_{t_j}
\label{Dhat}
\right)^2,
\end{eqnarray}
where the set $M_k = \{j: \one(X_{t_j},k) = 1\}$, and $card(M_k)=M$. 
Set $M_k$ is a set of indexes such that $X_{t_j} \in Bin_k$ and contains exactly $M$ time-instances.
The indicator function $\one(X_{t_j},k)$ is defined as
\[
\one(X_t,k) = \begin{cases}
  1, & X_{t} \in Bin_k,\\
  0, & X_{t} \not\in Bin_k,
  \end{cases}
\]
where 
$Bin_k = \left[x_k - {\dx}/{2}, x_k + {\dx}/{2} \right]$. 
Here, the indicator function  $\one(X_t,k)$ is analogous to conditioning in expressions \eqref{Ax}, \eqref{Bx}, but the conditioning is done on the interval $Bin_k$ instead of a particular value. We also impose that the $card(M_k) = M$
for all $k$, which means that we consider the situation when the number of time-instances for estimating the 
drift and diffusion coefficients does not depend on $x_k$. This implies that for all bins data always contains at least 
$M$ time-instances $t_j$ such that $X_{t_j} \in Bin_k$ for all $k$. In practice, such situation is likely to occur
when none of the $x_k$ are in the tails of the stationary distribution $\rho(x)$, e.g. $\max(x_k) - \min(x_k) \approx stddev(\rho(x))$. This is exactly the situation for many practical applications when the observational data is produced by numerical simulations or observations since rare events are unlikely to be a part of the trajectory $\{U_k, \, k=1,\ldots,N\}$.

In this paper we study analytical properties of estimations defined in \eqref{Ahat} and \eqref{Dhat}.
These estimators depend on three parameters - 
(i) the observational time-step $\dt$, (ii) the space-discretization $\dx$, and 
(iii) the number of observational time-instances $M$. Therefore, the key question is how to select these parameters to achieve optimal performance of estimators in \eqref{Ahat} and \eqref{Dhat}
while reducing the computational and data-generating complexities.  One obvious choice for selecting the parameters would be $\dt \to 0$, $\dx \to 0$, $M \to \infty$, but both $\dx \to 0$ and $M \to \infty$ increase the computational complexity of the problem. Moreover, if the observed data is fixed in size, it is not possible to achieve  $\dx \to 0$ and $M \to \infty$
simultaneously because as the width of the interval $Bin_k$ decreases, fewer observational points will satisfy 
$X_{t_j} \in Bin_k$. Therefore, in this paper we study the balance between three parameters, 
$\dt$, $\dx$, and $M$, which allows
achieving the optimal behavior of estimators in \eqref{Ahat} and \eqref{Dhat} with respect to the bias and the mean squared error.
To this end, we analyze the behavior of estimators as $\dt,~\dx \to 0$ and derive practical relationships
between $\dt$, $\dx$, and $M$ for small, but finite $\dt$ and $\dx$, in order to achieve optimal speed of convergence.

\subsection{Expectation with respect to the Truncated Density}

From the construction of the drift and diffusion estimators in \eqref{Ahat} and \eqref{Dhat} 
terms in the summation in the right-hand side of $\hat{A}(x_k)$ and $\hat{D}^2(x_k)$ are restricted 
to $X_{t_j} \in Bin_k$. 
Therefore, values of the stochastic process $X_{t_j}$ are sampled from a
stationary trajectory restricted to $Bin_k$.
Thus, we need to understand the stationary distribution restricted to $Bin_k$.
Formally such density can be represented as
\begin{equation}
\label{pk}
p_k(x) = G_k^{-1} \rho(x) \one(x,k) = 
\begin{cases}
G_k^{-1} \rho(x) & x \in Bin_k \\
0 & \text{otherwise},
\end{cases}
\end{equation}
where $G_k$ is the normalization factor
\[
G_k = \int\limits_{x_k - \dx/2}^{x_k + \dx/2} \rho(x) dx.
\]

To analyze the behavior of estimators $\hat{A}(x_k)$ and $\hat{D}^2(x_k)$ we first need to 
understand the asymptotic behavior (as $\dx \to 0$) 
of expectations with respect to the truncated density $p_k(x)$. 
For any function $f$ the expectation with respect to $p_k(x)$ is given by
\[
\bE_{p_k} [f(x)] = 
G_k^{-1} \int\limits_{x_k - \dx/2}^{x_k + \dx/2} f(x) \rho(x) dx.
\]
Considering sufficiently smooth functions $f$ and 
using Taylor expansions for $\rho(x)$ and $f(x)$ we obtain
\begin{equation}
\label{Epk}
\bE_{p_k} [f(x)] =
f(x_k) + \left[\frac{2f'(x_k)\rho'(x_k) + f''(x_k)\rho(x_k)}{\rho(x_k)}  \right] \dx^2/24 + O(\dx^4),
\end{equation}
which demonstrates explicitly the leading-order behavior of $\bE_{p_k} [f(x)]$.

\section{It\^o-Taylor Expansions}
\label{sec:sec3}

We utilize It\^o-Taylor expansions (see e.g. \cite{Peter})
to analyze the behavior of estimators \eqref{Ahat} and \eqref{Dhat} as $\Delta t \to 0$.
Assuming that $A(x)$ and $D(x)$ are sufficiently smooth functions, first few terms of
the It\^o-Taylor expansion of $X_{t_j + \dt}$ around $X_{t_j}$ can be written as 
\begin{eqnarray}
X_{t_j + \dt} &\approx& X_{t_j} + A(X_{t_j}) I_{(0), j} + D(X_{t_j}) I_{(1), j} +  B_2(X_{t_j}) I_{(1, 1), j} +
\nonumber \\
&&   B_3(X_{t_j}) I_{(0, 1), j} +  B_4(X_{t_j}) I_{(1, 0), j} +   B_5(X_{t_j}) I_{(0, 0), j} + B_6(X_{t_j}) I_{(1, 1, 1), j} =
\nonumber \\
&& X_{t_j} +  \sum_{q=0}^6 B_q(X_{t_j}) I_{\alpha_q, j} \equiv  ITE_j,
\label{ITE}
\end{eqnarray}
where we denote $B_0(x)\equiv A(x)$ and $B_1(x) \equiv D(x)$,
and other functions $B_k(x)$ are expressed through the drift and diffusion coefficients
\begin{align*}
& B_2(x) = D(x)D'(x), \quad B_3(x) = A(x)D'(x) + \frac{1}{2}D^2(x)D''(x), \\
& B_4(x) = D(x)A'(x), \quad
B_5(x) = A(x)A'(x) + \frac{1}{2}D^2(x)A''(x), \\
& B_6(X_{t_j}) = D(X_{t_j}) \left((D'(X_{t_j}) + D(X_{t_j})D''(X_{t_j})\right), 
\end{align*}
and $I_{\alpha_q,j}$ are stochastic integrals which are represented using indexes
$\alpha_q$
\begin{eqnarray*}
I_{(0), j} &=&  \int^{t_j + \dt}_{t_j} dt' = \dt, \qquad 
I_{(1), j} = \int^{t_j + \dt}_{t_j} dW_{t'},\\
I_{(0, 0), j} &=& \int^{t_j + \dt}_{t_j} \int^{s}_{t_j} dt' ds = \frac{\dt^2}{2}, \qquad 
I_{(0, 1), j} = \int^{t_j + \dt}_{t_j} \int^{s}_{t_j} dt' dW_s ,\\
I_{(1, 0), j} &=&  \int^{t_j + \dt}_{t_j} \int^{s}_{t_j} dW_{t'} ds, \qquad 
I_{(1, 1), j} =  \int^{t_j + \dt}_{t_j} \int^{s}_{t_j} dW_{t'} dW_s,\\
I_{(1, 1, 1), j} &=& \int^{t + \Delta t}_{t} \int^{s}_{t} \int^{t'}_{t} dW(r) dW(t') dW(s).
\end{eqnarray*}
Index $\alpha_q$ determines the order of integration in stochastic integrals.
1 in index $\alpha_q$ corresponds to integration with respect to the Brownian motion, and 0 corresponds to integration with respect to time.
Therefore, from definition \eqref{ITE}, 
$\alpha_0 = (0)$, $\alpha_1 = (1)$, $\alpha_2 = (1,1)$, etc.
Properties of these stochastic integrals have been studied, for example, in \cite{Peter}.
Integrals $I_{(0), j}$ and $I_{(0, 0), j}$ are deterministic, while others are random variables.
Integrals $I_{(1), j}$, $I_{(0, 1), j}$, and $I_{(1, 0), j}$ are Gaussian with mean zero and variances
\begin{equation}
\label{L2I1}
\bE [I_{(1), j}^2]=\dt, \quad \bE[I_{(0, 1), j}^2]=\dt^3/3, \quad \bE[I_{(1, 0), j}^2]=\dt^3/3.
\end{equation}
Integrals $I_{(1, 1), j}$ and $I_{(1, 1, 1), j}$ are non-Gaussian with the first two moments given by
\begin{eqnarray}
\label{L2I2}
&& \bE\left[I_{(1, 1), j}\right] =0, \qquad \bE[I_{(1, 1), j}^2] = \dt^2/2, 
\label{L2I2} \\
&& \bE\left[I_{(1, 1, 1), j}\right] =0, \qquad \bE[I_{(1, 1, 1), j}^2] = O(\dt^3).
\label{L2I3}
\end{eqnarray}
Moreover, one can prove that
\begin{equation}
    I_{(1,1),j} = \left( (\Delta W_{j+1})^2 - \dt \right)/2,
    \label{I11}
\end{equation}
where  $\Delta W_{j+1} = W_{t_j + \dt} - W_{t_j}$. 
Mixed second moments of stochastic integrals are
\begin{align}
& \bE[I_{(1), j} I_{(1,0), j}] \le \frac{\dt^2}{2} , \quad 
\bE[I_{(1), j} I_{(0,1), j}]  \le \frac{\dt^2}{2}, \quad
\bE[I_{(1,0), j} I_{(0,1), j}] \le \frac{\dt^2}{2},
\nonumber \\
& \bE[I_{(1), j} I_{(1,1), j}] =
\bE[I_{(1,1), j} I_{(0,1), j}]  =
\bE[I_{(1,1), j} I_{(1,0), j}] = 0. 
\label{Ib1} \\
& \bE[I_{1,1}^4,j] = O(\dt^4).
\nonumber 
\end{align}
Triple stochastic integrals with even number of ones are of higher order and do not make any low-order contributions in calculations of the bias and the mean-squared error discussed in subsequent sections.

Following \cite{Peter} it is also useful to introduce function $n(\alpha_q)$ and 
$n(\alpha_q, \alpha_l)$ which counts the number of ones
\begin{align}
\label{n}
\begin{split}
n(\alpha_q) &= \text{number of ones in~} \alpha_q \\
n(\alpha_q,\alpha_l) &= \text{number of ones in~} \alpha_q \text{~and~} \alpha_l.
\end{split}
\end{align}

\section{Bias of the Drift and Diffusion Estimators}
\label{sec:sec4}
In this section 
we analyze the bias for the drift and diffusion estimators in \eqref{Ahat}, \eqref{Dhat}.
We show that these estimators are biased for finite $\dt >0$ and $\dx > 0$, but the bias vanishes 
in the limit $\dt\to 0$ and $\dx\to 0$.

\subsection{Bias of $\hat{A}(x_k)$}
\label{sec:biasA}

To analyze the bias of $\hat{A}(x_k)$ we consider the expected value of $\hat{A}(x_k)$
\begin{align*}
\bE[\hat{A}(x_k)] = & \frac{1}{M\dt}    
\sum\limits_{j \in M_k} \bE \left[X_{t_j + \dt} - X_{t_j}|X_{t_j}\in Bin_k\right] \approx \\
& \frac{1}{M\dt}  \sum\limits_{j \in M_k} \bE \left[ITE_j - X_{t_j}|X_{t_j}\in Bin_k\right]  = \\
& \frac{1}{M\dt}   \sum\limits_{j \in M_k} 
\sum\limits_{l=0}^6\bE [B_l(X_{t_j}) I_{\alpha_l, j}|X_{t_j}\in Bin_k].
\end{align*}

If we denote the filtration generated by $W_t$ as $\calF_t$ then
\begin{align}
\label{expect}
& \bE[B_l(X_{t_j}) I_{\alpha_l, j}|X_{t_j}\in Bin_k] = 
\bE\left[ \bE[B_l(X_{t_j}) I_{\alpha_l, j}|\calF_{t_j}] | X_{t_j}\in Bin_k\right] = 
\nonumber \\
& \bE \left[ B_l(X_{t_j}) | X_{t_j}\in Bin_k \right] \bE[I_{\alpha_l, j}] = 
\bE_{p_k} \left[ B_l(x)\right] \bE[I_{\alpha_l, j}]
\nonumber
\end{align}
and we can use properties of stochastic integrals to evaluate $\bE[I_{\alpha_l, j}]$.
We would like to point out that conditional expectation 
$\bE_{p_k} \left[ B_l(x)\right]$ in general depends on $x_k$.
Thus, we obtain
\begin{align*}
\bE[\hat{A}(x_k)] \approx & \frac{1}{M\dt} \sum\limits_{j \in M_k} \Big[ \bE_{p_k}[B_0(x)] \dt +
\bE_{p_k} [B_2(x)] \dt^2/2 + \bE_{p_k} [B_5(x)] \dt^2/2  \Big] = \\
& A(x_k) +  O(\dx^2) + O(\dt).
\end{align*}
%
Therefore, 
\[
\bE[\hat{A}(x_k)]  \to A(x_k) \text{~~as~~} \dt, \, \dx \to 0.
\]

For small, but finite $\dx$ and $\dt$, we can expect that 
\begin{equation}
\label{bias}
Bias[\hat{A}(x_k)] \sim C(\dx^2 + \dt),
\end{equation}
where  constant $C \equiv C(x_k)$ might depend on $x_k$. Therefore, formula \eqref{bias}
indicates that in order to balance the bias terms on the right-hand side of \eqref{bias}
the space- and time-discretization should scale as
\[
\dx^2 \sim \dt.
\] 
This scaling has important practical implications indicating that the bin size can be taken to be quite large
compared to the observational time-step, $\dt$. We will discuss the scaling between $\dx$ and $\dt$ further in other sections.

\subsection{Bias of $\hat{D}^2(x_k)$}

We analyze bias of $\hat{D}^2(x_k)$ in a manner similar to the previous section.
We consider
\begin{align*}
\bE[\hat{D}^2(x_k)] = & \frac{1}{M\dt}    
\sum\limits_{j \in M_k} \bE \left[ (X_{t_j + \dt} - X_{t_j})^2 | X_{t_j} \in Bin_k\right] \approx \\
& \frac{1}{M\dt}    
\sum\limits_{j \in M_k} \bE \left[ (ITE_j - X_{t_j})^2 | X_{t_j} \in Bin_k\right] = \\
& \frac{1}{M\dt}  \sum\limits_{j \in M_k} 
\sum\limits_{l, q=0}^6 \bE [B_l(X_{t_j}) B_q(X_{t_j}) I_{\alpha_l, j} I_{\alpha_q, j} |X_{t_j}\in Bin_k] = \\
& \frac{1}{M\dt}  \sum\limits_{j \in M_k} 
\sum\limits_{l, q=0}^6 \bE [B_l(X_{t_j}) B_q(X_{t_j})  |X_{t_j}\in Bin_k] \, \bE[I_{\alpha_l, j} I_{\alpha_q, j}].
\end{align*}
Thus, we need to compute expected values of cross-products $\bE[I_{\alpha_l, j} I_{\alpha_q, j}]$ for all
$l,q = 0,\ldots,6$. Deterministic terms resulting from $I^2_{(0),j}$ and $I^2_{(0,0),j}$ are non-zero, but 
$I^2_{(0,0),j}$ is of higher order.
For stochastic terms, we can show that terms where $n(\alpha_q, \alpha_l)$ is odd are zero.
Thus, the leading-order term arises from $\bE[I_{(1), j}^2] = \bE[(W_{t_j + \dt} - W_{t_j})^2]=\dt$. 
Other non-zero terms (see \eqref{L2I1}, \eqref{L2I2}, \eqref{Ib1}) are of higher order.
Therefore, we obtain,
\begin{align*}
\bE[\hat{D}^2(x_k)] = \frac{1}{M\dt}    
\sum\limits_{j \in M_k} \bE_{p_k} [D^2(x)] \dt + O(\dt^2) = 
D^2(x_k) + O(\dx^2) + O(\dt).
\end{align*}
Similar to the drift estimator, 
$Bias[\hat{D}^2(x_k)] \sim C(\dx^2 + \dt)$ and
diffusion estimator becomes unbiased in the
limit $\dt, \dx \to 0$.
However, we would like to emphasize that the above scaling is applicable to the diffusion squared, 
$D^2(x_k)$, not $D(x_k)$.

\subsection{Comments on another possible drift estimator}

One can define a slightly different drift estimator (c.f. with $\hat{A}(x_k)$ in \eqref{Ahat})
\begin{equation}
\label{Ahhat}
{\tilde{A}} (x_k) = 
\frac{1}{M\dt} \sum_{j \in M_k} (X_{t_j + \dt} - x_k).
\end{equation}
In this case, the estimator is centered at $x_k$, instead of subtracting $X_{t_j}$.

We can compute the bias of the drift estimator in \eqref{Ahhat} in a manner totally similar to the computations
of the bias for $\hat{A}(x_k)$ in section \ref{sec:biasA}, i.e.
\begin{align*}
\bE[{\tilde{A}}(x_k)] = & \frac{1}{M\dt}  
\sum\limits_{j \in M_k} \bE \left[X_{t_j + \dt} - x_{k}|X_{t_j} \in Bin_k\right] = \\
& \frac{1}{M\dt} \sum\limits_{j \in M_k} \bE \left[ITE_j - x_{k}|X_{t_j} \in Bin_k\right] + O(\dt) = \\
& \bE_{p_k}[\hat{A}(x_k)] + \frac{1}{\dt} \frac{1}{M}    \sum\limits_{j \in M_k} \bE_{p_k}[x - x_k] + O(\dt) = \\
& \bE_{p_k}[\hat{A}(x_k)] + r(x_k) \frac{\dx^2}{\dt} + O(\dt) = A(x_k) + r(x_k) \frac{\dx^2}{\dt} + O(\dx^2) + O(\dt),
\end{align*}
where we can compute the remainder by applying \eqref{Epk} with 
$f(x) = x - x_k$ and, therefore,
$r(x_k) = \rho'(x_k)/(12 \rho(x_k))$.
Thus, there is an additional condition $\dx^2/\dt \to 0$ for the estimator 
in \eqref{Ahhat} to be asymptotically unbiased.
In addition, the term $r(x_k) \dx^2/\dt$ can provide a significant contribution to the bias of the estimator
\eqref{Ahhat} for finite $\dx$ and $\dt$. Similar issue arises if we consider modified estimator for the diffusion coefficient.
Thus, estimator \eqref{Ahhat} is inferior compared to the estimator \eqref{Ahat} and we will consider
\eqref{Ahat} for the rest of this paper.

\section{MSE of the Drift and Diffusion Estimators}
\label{sec:sec5}

Next, we compute the leading-order behavior of the Mean-Squared-Error (MSE) for both estimators $\hat{A}(x_k)$ and $\hat{D}^2(x_k)$.
The calculation is quite technical, especially for the diffusion estimator, 
and we only sketch here most important points. Details are presented in the Appendix.

\subsection{MSE of the Drift Estimator}
\label{sec:msedrift}
In order to understand the behavior of the Mean-Squared-Error for the drift estimator we need to compute 
$\|\hat{A}(x_k) - A(x_k)\|_2^2$
where the norm is computed conditioned on $X_{t_j} \in Bin_k$ for all $j\in M_k$.
In particular, we compute
\begin{align}
\| \hat{A}(x_k) & - A(x_k) \|_2^2 = \bE \left[ \left(\hat{A}(x_k) - A(x_k) \right)^2 \Big| X_{t_j} \in Bin_k\right]= 
\label{MSEADEF} \\
& \bE \left[ \left(\frac{1}{M\dt}  \sum\limits_{j \in M_k} (X_{t_j + \dt} - X_{t_j}) - A(x_k) \right)^2 \Big| X_{t_j}\in Bin_k \right] \approx
\nonumber \\
& \bE \left[ \left(\frac{1}{M\dt}  \sum\limits_{j \in M_k} \left( 
\left[A(X_{t_j}) - A(x_k) \right]\dt + \sum\limits_{q=1}^6 B_q(X_{t_j})I_{\alpha_q, j} \right)\right)^2 \Big| X_{t_j}\in Bin_k \right].
\label{MSEA}
\end{align}
The details of calculating the expectation in \eqref{MSEA} are presented in Appendix \ref{ap1}.

For the rest of the paper we will use ``$C$'' to denote a generic constant.
Since our analysis is local (i.e. restricted to a particular $x_k$ and $Bin_k$)
this constant might depend on $x_k$ and $\dx$, but for each bin this constant converges to a finite 
value as the bin size goes to zero, i.e.
$C \to C_{lim}(x_k) > 0$ as $\dx \to 0$.

We obtain the following asymptotic result for the MSE of the drift estimator
\begin{equation}
\label{MSEA1}
\| \hat{A}(x_k) - A(x_k) \|_2^2 \le C \left( \sqrt{\dt} + \frac{1}{M\dt} + (\dx)^2 + \frac{\dx}{\sqrt{\dt}} \right) + 
h.o.t.,
\end{equation}
where higher-order-terms involve various higher-order powers of $\dt$ and $\dx$.
First, we notice that the requirements for the $MSE\{\hat{A}(x_k)\} \to 0$ are
\begin{equation}
\label{cond1}
M\dt \to \infty, \quad
\dt \to 0, \quad
\dx \to 0, \quad
\frac{\dx}{\sqrt{\dt}} \to 0.
\end{equation}
The first three conditions are expected, but the last condition provides a relationship between $\dt$ and $\dx$. 
Asymptotic behavior of the $MSE\{\hat{A}(x_k)\}$ confirms the optimal relationship between
$\dx$ and $\dt$ derived for the bias of 
$\hat{A}(x_k)$ in section \ref{sec:biasA}. 
To guarantee $MSE\{\hat{A}(x_k)\} \to 0$,
the spatial discretization can be chosen to be much coarser
than the observational time-step and the appropriate practical scaling is
\begin{equation}
\label{dxdtscale_diff}
\dx \sim \dt^{1/2 + \eps},
\end{equation}
where $\eps$ is any fixed small number.
This scaling motivated by the fact that in practice we would like the bin size to be as large as possible.
Larger bin sizes allow increasing the number of points which fall in each bin and thus $card(M_k)$ becomes larger.

\subsection{MSE of the Diffusion Estimator}
\label{sec:msediff}
To consider the MSE of the diffusion estimator we need to compute the leading-order behavior of 
the following conditional expectation
\begin{align}
\| \hat{D}^2(x_k) & - D^2(x_k) \|_2^2 = \bE \left[ \left(\hat{D}^2(x_k) - D^2(x_k) \right)^2 \Big| X_{t_j} \in Bin_k\right]= 
\label{MSEDDEF} \\
& \bE \left[ \left(\frac{1}{M\dt}  \sum\limits_{j \in M_k} (X_{t_j + \dt} - X_{t_j})^2 - D^2(x_k) \right)^2 \Big| X_{t_j}\in Bin_k \right] \approx
\nonumber \\
& \bE \left[ \left(\frac{1}{M\dt}  \sum\limits_{j \in M_k} \left(
\left( \sum\limits_{l=0}^6 B_l I_{\alpha_l,j} \right)^2 
- D^2(x_k) \dt \right) \right)^2 \Big| X_{t_j}\in Bin_k \right].
\label{MSED}
\end{align}
This expression contains many terms, but most of them can be treated in a similar manner.
Moreover, 
to obtain asymptotic behavior of the $MSE$ for the diffusion estimator we only need 
to keep track of lowest-order terms which typically arise from 
the first few stochastic integrals. Details of the calculation are presented in Appendix \ref{ap2}.
The asymptotic behavior of the MSE for the diffusion estimator is given by 
\begin{equation}
\label{MSED1}
\| \hat{D}^2(x_k) - D^2(x_k) \|_2^2 = C \left(\frac{1}{M} +\dx + \dt \right)+ h.o.t.
\end{equation}
The asymptotic behavior of the MSE for the diffusion estimator is different from the 
MSE of the drift estimator which is consistent with results for other estimators, 
such as the Maximum Likelihood Estimators \cite{degia2012} which exhibit different convergence rates.
Asymptotic conditions for the $MSE\{\hat{D}^2(x_k)\}\to 0$ are less demanding than for the drift estimator. In particular, the diffusion can be accurately estimated on finite-time intervals $M\dt = T = Const < \infty$. 
The optimal relationship between the space- and time-discretizations is 
\begin{equation}
\label{dxdtscale}
\dx \sim\dt,
\end{equation}
which is different compared to the optimal scaling for the drift estimator. 
However, it is difficult to access analytically the value of constants multiplying $\dx$ and $\dt$ in the expression 
for the MSE of the diffusion estimator. Magnitudes of these constants might play an important role in practice
and might lead to a different scaling regime for practical values of $\dt$ under consideration. We address this
issue numerically in the next section.

\section{Numerical Simulations}
\label{sec:sec6}
In this section we perform numerical simulations and analyze numerically the validity of 
expressions for the MSE of the drift and diffusion estimations in 
\eqref{MSEA1} and \eqref{MSED1}, respectively. 
%
%
In particular, MSE estimates for the drift and diffusion in \eqref{MSEA1} and \eqref{MSED1} 
indicate very different behavior when $M\to\infty$ with respect to two asymptotic regimes 
$M \dt = Const$ and $M \Delta t \to \infty$. 
Moreover, our analytical results also indicate that we can use quite large 
$\dx$ for adequate estimation of the drift. We will show numerically that
larger bin size, $\dx$, can be used in the estimation of both, the drift and the diffusion coefficients.
We also investigate the role of the observational time-step, $\dt$, which appears in the
denominator in \eqref{MSEA1}. Finally, we demonstrate that nonlinear regression 
can recover the correct form of the drift and diffusion coefficients.

\subsubsection*{Numerical Simulations for the Cubic Process}
We illustrate behavior of estimators with respect to changing 
$M$, $\dt$, and $\dx$ using numerical data 
generated by the following SDE with cubic drift and linear diffusion
\begin{equation}
\label{eq:cubic}
dX_t = -\gamma X_t^3 dt + \left( \sigma_1 + \sigma_2 X_t \right) dW_t
\end{equation}
with 
parameters $\gamma = 1$ and $\sigma_1= \sigma_2=1/\sqrt{2}$. 

In the first regime $M\dt = Const$ we use the following values
\begin{equation}
\label{mdtconst}
(M,\dt) = 
(50, 0.02), \, (100, 0.01), \, (200,0.005), \, (500,0.002), \, (1000,0.001)
\end{equation}
and in the second regime $M\dt \to \infty$ we use 
\begin{equation}
\label{mdtinf}
\dt=0.01 \quad \text{and} \quad 
M=50, \, 100, \, 200, \, 500, \, 1000.
\end{equation}
The drift and diffusion coefficients are estimated on a discrete mesh 
$x_k \in [-L,L]$ with $L=0.5$ and 
the Number of Bins is $NB=10, 20, 40$ which corresponds to $\dx = 0.1$, $0.05$, $0.025$.
The choice of $L=0.5$ is motivated by the fact that the stationary standard deviation 
of the process in \eqref{eq:cubic} with chosen parameters is approximately 0.5.
One can choose to estimate the drift and diffusion coefficients on a larger interval,
but the data for large $L$ would become scarce, since $|x_k| \gg 0.5$ corresponds to
values of $X_t$ in the tails of the stationary distribution.
Thus, a near-optimal practical guideline is to estimate the drift and diffusion coefficients on an interval $Mean \mp StdDev$.
We use the 1.5 strong discretization (see \cite{Peter}) to generate stationary trajectories of the SDE in \eqref{eq:cubic}. The time-step of integration is $\delta t=5 \times 10^{-4}$.

To compute the MSE numerically we perform Monte-Carlo simulations and compute many realizations of
sampled trajectories and, in turn, of the drift and diffusion estimators. 
Then, we compute the discrete analog of the MSE
 \begin{eqnarray}
MSE_{drift} &=& \frac{1}{MC}\sum^{MC}_{j = 1}\left( \sum_k \left(\hat{A_k}^{(j)} - A(k)\right)^2\Delta x\right), \nonumber \\
MSE_{diff} &=& \frac{1}{MC}\sum^{MC}_{j = 1} \left( \sum_k \left(\hat{D^2_k}^{(j)} - D^2(k)\right)^2\Delta x\right), \nonumber
\end{eqnarray}
where
$MC$ is the number of Monte-Carlo Realizations, $\hat{A_k}^{(j)}$ and $\hat{D^2_k}^{(j)}$ 
are the drift and diffusion estimators computed for the $j$-th Monte-Carlo realization, 
and $k$ represents the $k$-th bin (all bins are of size $\dx$). We use $MC=500$ trajectories.

\begin{figure}[ht!]
    \centerline{
    \includegraphics[scale=0.5]{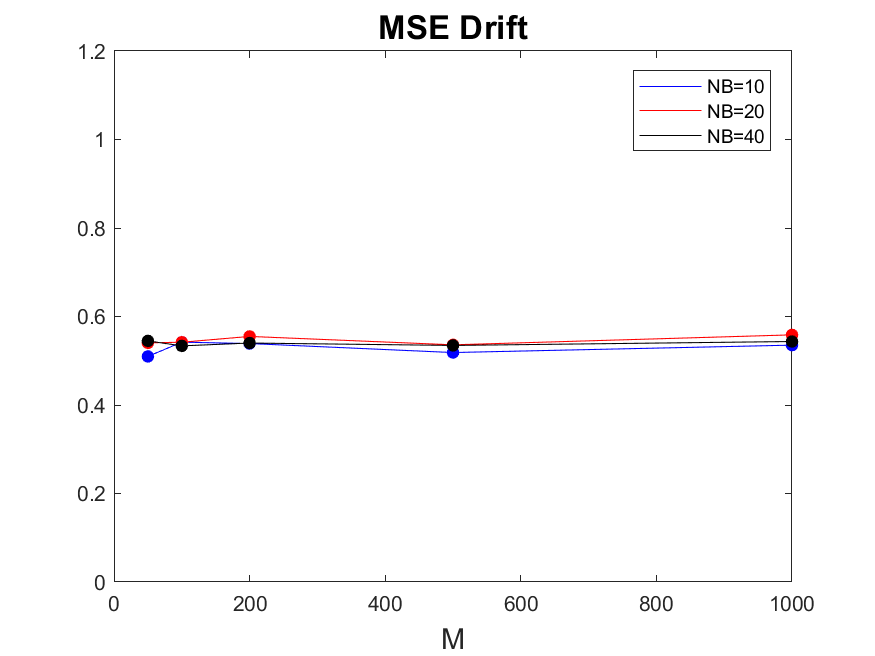} \hspace{-0.7cm}
    \includegraphics[scale=0.5]{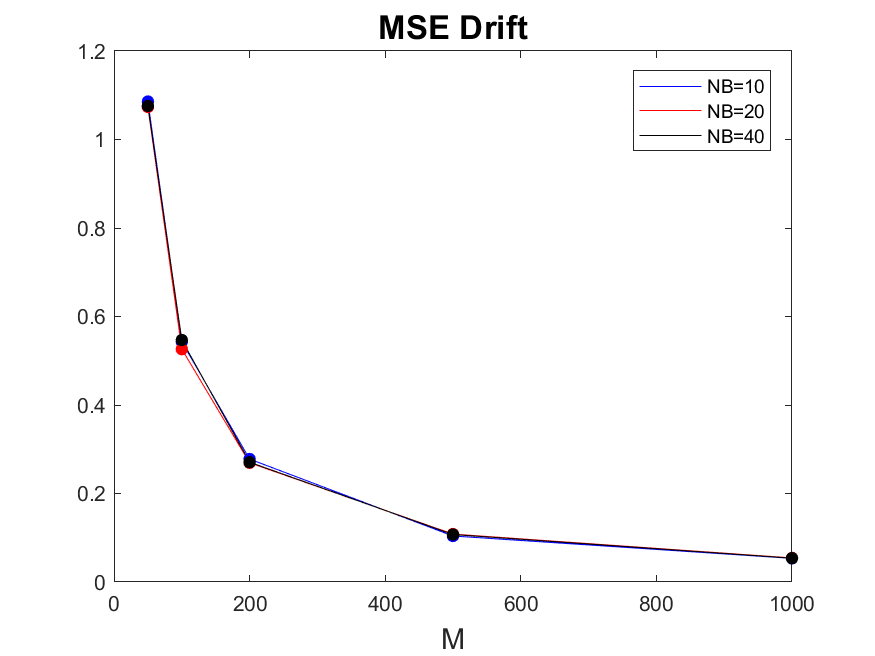}}
    \centerline{
    \includegraphics[scale=0.5]{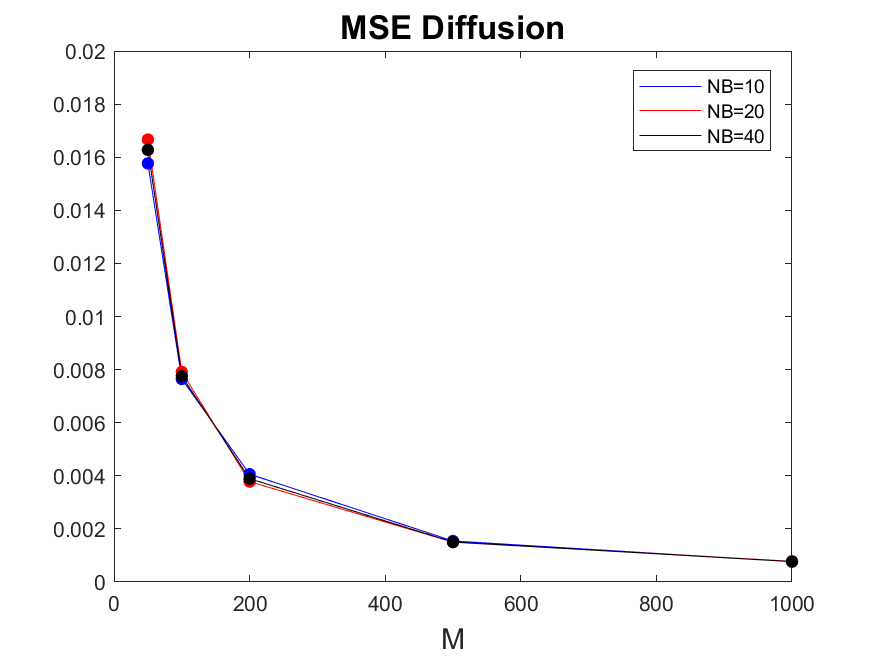} \hspace{-0.7cm}
    \includegraphics[scale=0.5]{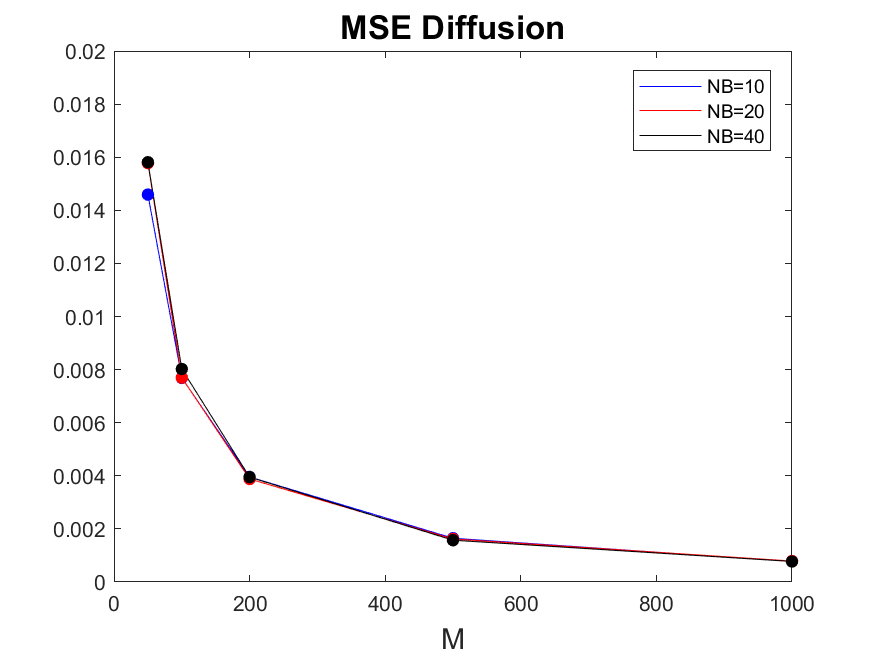}}
    \caption{MSE of the Drift (top)  and Diffusion (bottom) estimators 
    with two different sampling regimes 
    $M\dt = 500$ (left) and $M\dt \rightarrow \infty$ (right).
    Corresponding parameters are given by \eqref{mdtconst} and \eqref{mdtinf}
    for $M\dt = 500$ and $M\dt \rightarrow \infty$, respectively.
    Simulations of the cubic process \eqref{eq:cubic} with 
    $\Delta x = 1/NB$. Three values of $\Delta x=0.1$, $0.05$, $0.025$ 
    overlap almost completely on all plots.
   }
\label{fig:MSE}
\end{figure}

Figure \ref{fig:MSE} depicts the behavior of the $MSE$ for the drift and diffusion coefficients.
There is a clear evidence that numerical errors for the drift and diffusion coefficients 
behave very differently for the two sampling regimes $M\dt \to \infty$ and $M\dt = Const$. 
In particular, errors for the diffusion 
estimator are decaying as long as $M\to\infty$. 
However, the behavior of the drift estimator depends drastically
on whether $M\Delta t= Const$ or $M\dt \to\infty$. When $M\dt= Const$, the MSE for the drift estimator remains constant as $M\to\infty$ (and $\dt\to 0$) as predicted by our analytical expression in \eqref{MSEA1}. When $M\dt \to\infty$ we observe decay of errors for
the drift estimator as $M\to\infty$ (and $\dt$ fixed).
Moreover, the slope of the $MSE$ vs $M$ on the log-log plots (not depicted here for 
the brevity of presentation) for the diffusion estimator equals approximately $-1$.

\begin{figure}[ht!]
\centerline{\includegraphics[scale=0.5]{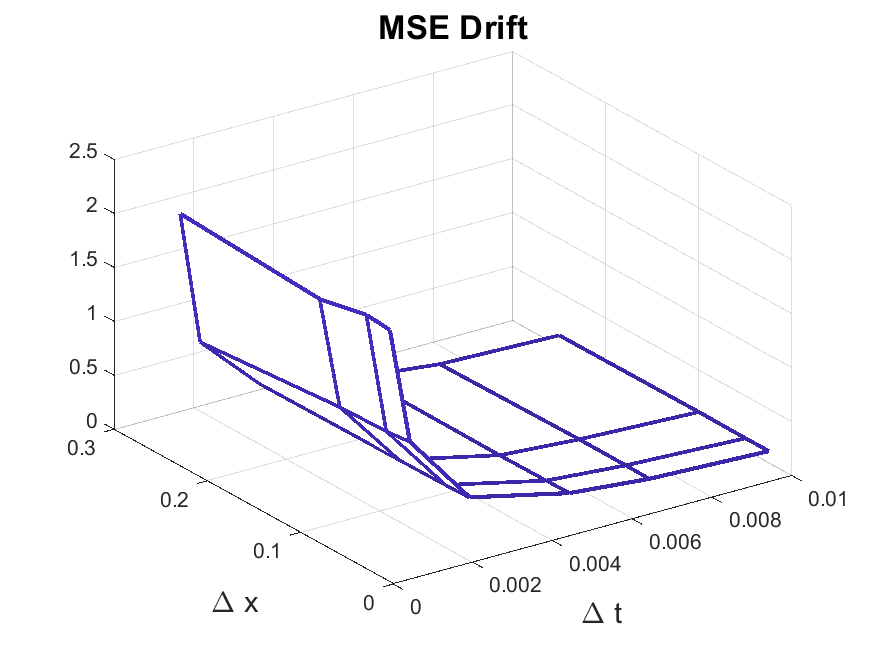}
\includegraphics[scale=0.5]{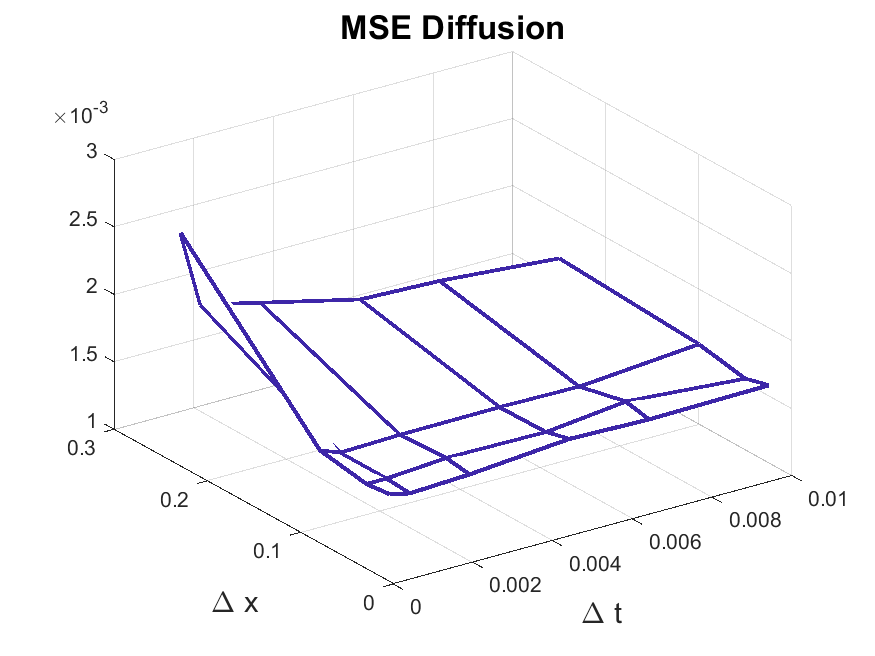}}
\caption{Averaged (over all bins) Mean-Squared-Errors for the estimation of the drift (left) and diffusion (right) coefficients. Simulations of the cubic process \eqref{eq:cubic} with $M=500$ and $MC=500$.}
\label{fig:l23d}
\end{figure}

\begin{figure}[ht!]
\centerline{\includegraphics[scale=0.5]{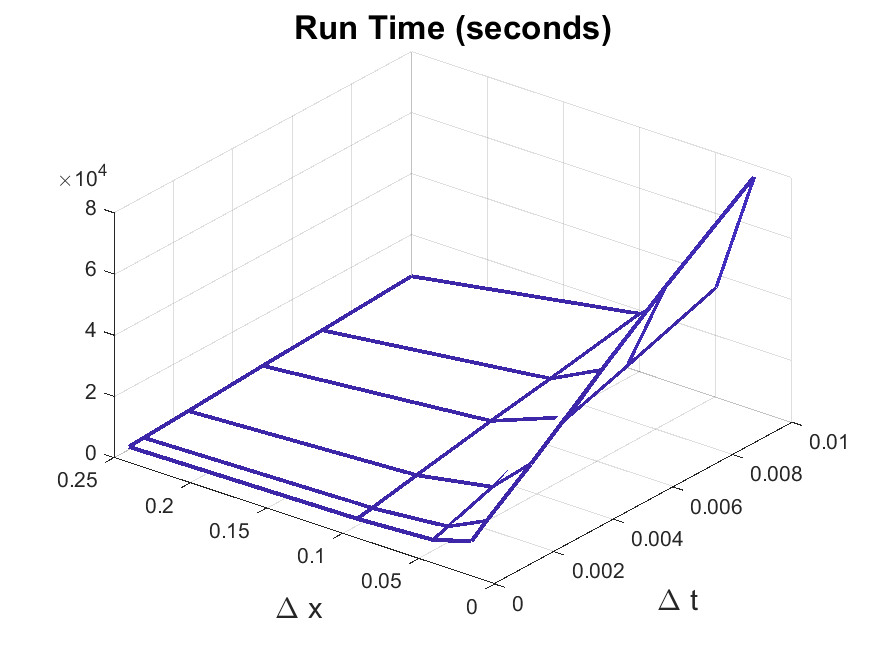}}
\caption{Running time (seconds) to generate data for estimation of drift and diffusion coefficients for the estimation of the drift and diffusion coefficients. Simulations of the cubic process \eqref{eq:cubic} with $M=500$ and $MC=500$.}
\label{fig:runt}
\end{figure}

It is somewhat difficult to disentangle contributions of errors from different terms outlined in our analytical expressions 
\eqref{MSEA1} and \eqref{MSED1}.
To analyze the behavior of the Mean-Squared-Errors in the $\dt - \dx$ plane, we performed simulations for a range of these parameters. 
In particular, we choose 
\[
\dt = 0.0005, 0.001, 0.0025, 0.005, 0.007, 0.01, \quad
\dx = 0.025, 0.05, 0.1, 0.25,
\]
but keep the number of points in each bin $M=500$.
We present averaged (over all bins) MSE errors in Figure \ref{fig:l23d}, we also present the running time required to generate the corresponding datasets in Figure \ref{fig:runt}. We would like to point out that we generate the data for estimation with non-overlapping sampling time-step $\dt$ and the code terminates when there are $M=500$ in each bin.
\begin{figure}[ht!]
\centerline{\includegraphics[scale=0.55]{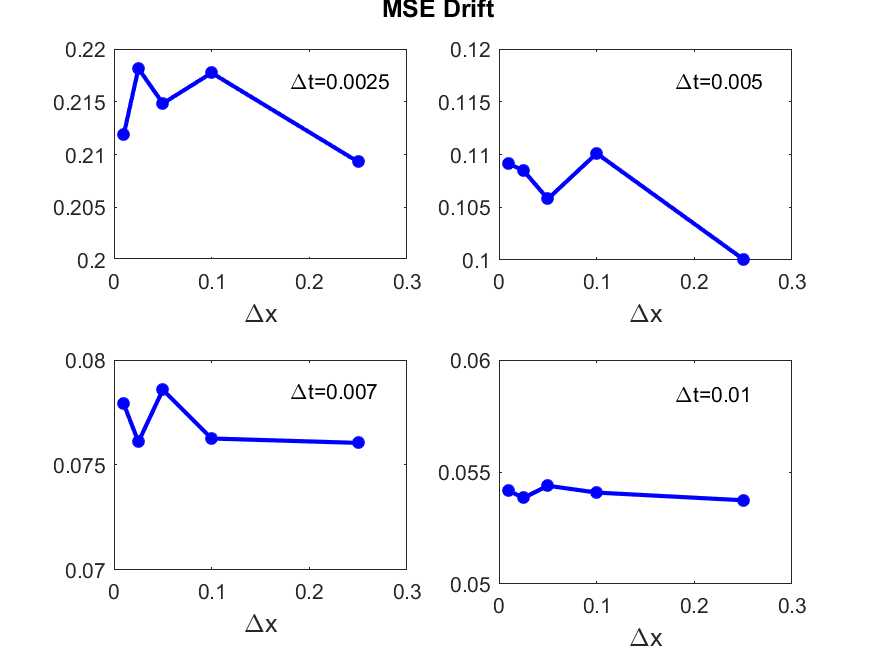}
\hspace*{-0.5cm}
\includegraphics[scale=0.55]{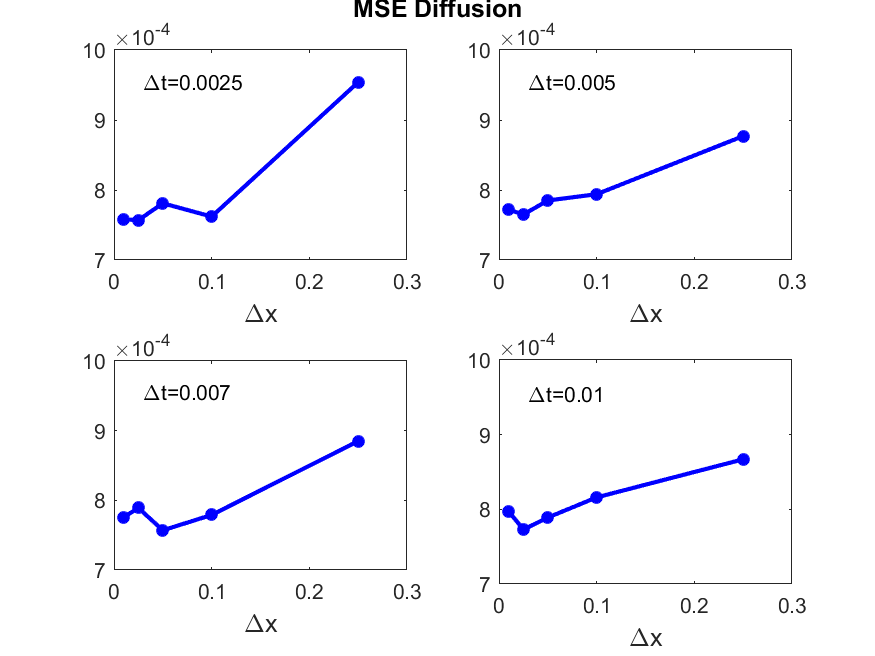}}
\caption{Averaged (over all bins) Mean-Squared-Errors for the estimation of the drift (left) and diffusion (right) coefficients vs $\Delta x$ for several particular values of $\Delta t=0.0025$, $0.005$, $0.007$, $0.01$. Simulations of the cubic process \eqref{eq:cubic} with $M=500$ and $MC=500$.}
\label{fig:l2snap}
\end{figure}
Figure \ref{fig:l23d} indicates that errors of estimation are not really affected by the choice of $\dx$. 
There is a slight error increase for $\dx=0.25$ for the diffusion estimator, but, errors in the estimation of the diffusion term 
remain very small for all $\dx$ and $\dt$.
Thus, there is an indication that the diffusion estimator is affected primarily by $\dt$ only
for $\dx=0.25$, with error increasing for smaller $\dt$.
However, this increase is rather small (about $1.6\times 10^{-3}$ for $\dt=0.01$ vs $2.6\times 10^{-3}$ for $\dt=0.0005$).
Overall, errors for the diffusion estimator remain approximately two orders of magnitude smaller compared to the errors 
for the drift estimator. Smallest errors for the drift estimator of order $O(10^{-1})$ are reached at $\dt=0.1$.
In addition, we also performed simulations for the same choice of parameters $\dt$ and $\dx$, but with $M=1000$.
In simulations with $M=1000$ all errors (both drift and diffusion estimation) become uniformly approximately twice smaller. This is
a strong indication that error terms which involve $(M\dt)^{-1}$ and $M^{-1}$ for the drift and diffusion estimations in 
 \eqref{MSEA1} and \eqref{MSED1}, respectively, are dominant and other terms in \eqref{MSEA1} and \eqref{MSED1} 
 do not contribute significantly in the parameter regime considered here.

From these simulations we can conclude the following - (i) drift estimator is affected significantly by the choice of $\dt$
and is not affected by $\dx$; since we do not observe significant variations of error for the drift estimator with respect to $\dx$,
the leading order error term for the drift estimator is $(M\dt)^{-1}$; 
(ii) diffusion estimator is not affected significantly by either $\dt$ or $\dx$; therefore, 
the leading error term for the diffusion estimator is $M^{-1}$; 
(iii) for practical values of $M=500,\ldots,1000$ other error terms in \eqref{MSEA1} and \eqref{MSED1} do not seem to be significant; 
(iv) numerical results presented here indicate that it is beneficial to select $\dx$ and $\dt$ to be
relatively large with $\dx \gg \dt$; 
(v) Figure \ref{fig:runt} also indicates that it is beneficial to select larger
$\dx$ to minimize computational or experimental data-generating effort since computational time increases significantly for small $\dx$.

We would like to point out that here estimation of the drift and diffusion is carried out on a
non-overlapping spatial grid. While it is possible to utilize overlapping bins and, thus,
use the same observational pair of points to perform estimation of the drift and diffusion for two neighboring bins, it is likely that this approach would result in correlated errors in two neighboring bins and, therefore, would not provide any additional benefits when using regression to estimate the drift and diffusion.

We also would like to point out that we consider a large range of bin sizes  $\dx \in [0.01, 0.25]$,
but numerical errors are almost completely insensitive to the 
variations of the bin size in this range. Therefore, the choice of $\dx$ should be motivated by the required number of points for performing
a regression fit after the non-parametric estimation. In most cases, 10 to 20 spatial 
points should be enough to preform regression estimation for non-oscillatory functions.
From Figure \ref{fig:l23d} we can see that drift and diffusion estimation with 
$\dx=0.1$ and $\dt=0.01$ provides adequate results.
Therefore, numerical results presented here suggest that the scaling in \eqref{dxdtscale_diff} 
should be appropriate in many practical situations.

\begin{figure}[ht!]
    \centerline{
   \includegraphics[scale=0.5]{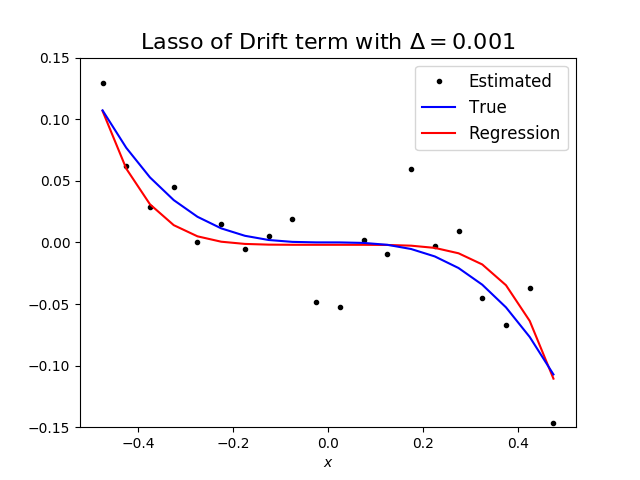} \hspace{-0.7cm}
    \includegraphics[scale=0.5]{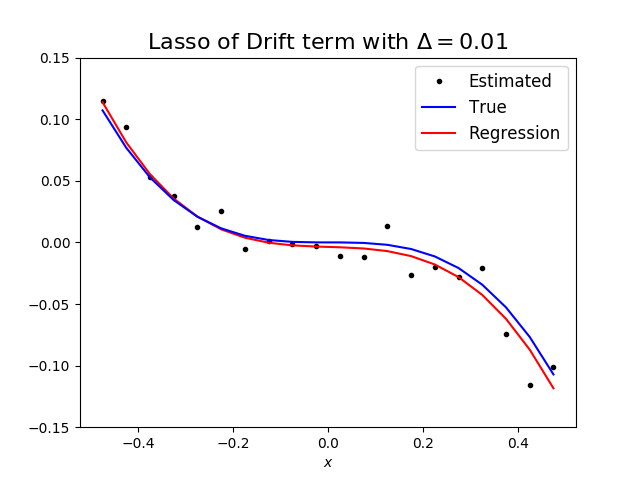}}
    \centerline{
   \includegraphics[scale=0.5]{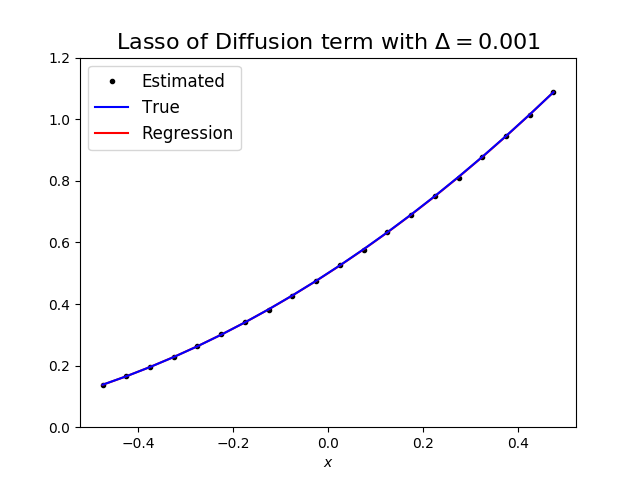} \hspace{-0.7cm}
    \includegraphics[scale=0.5]{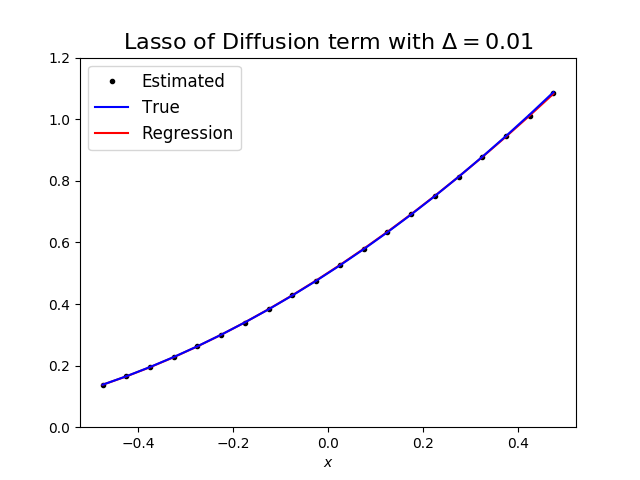}}
       \caption{Lasso regression fit for the drift (top part) and diffusion (bottom part) terms
       for the estimation of the cubic process \eqref{eq:cubic} with $\Delta x = 0.05$
  and $(M,\dt) = (1000, 0.001)$  (left)
   and $(M,\dt) = (1000, 0.01)$ (right). Coefficients of the fitted polynomials are presented in Table \ref{table:reg}.}
\label{fig:reg}
\end{figure}
Next, we perform nonlinear regression to recover the functional form of the drift and diffusion 
using discretely estimated data $\hat{A_k}$ and $\hat{D_k^2}$. 
We compare results of polynomial regression for two sets of parameters 
$(M,\dt)=(1000,0.001)$ and  $(M,\dt)=(1000,0.01)$ with $\dx=0.05$.
Since in practice we do not know
the functional form of the drift and diffusion in advance, we tested several regression
techniques for fitting high degree polynomials. Here we present results of fitting
the seventh degree polynomial. 
We also obtained similar results for fitting polynomials of degree 9.
In particular, we used standard nonlinear regression,
Lasso ($L^1$ penalty terms), and Ridge ($L^2$ penalty terms) regularizations.
Results are presented in Figure \ref{fig:reg} and Table \ref{table:reg}.
Lasso regression resulted in a better estimation of coefficients compared with
polynomial and Ridge regression.
In particular, standard polynomial 
regression is unstable, resulting in large coefficients for all powers.
Ridge regression is often similar to Lasso regression, but Lasso appears to be more stable resulting
in robust estimation for a wide range of the penalty parameter and 
polynomials of different order.
Here we depict only results for the Lasso regression in 
Figure \ref{fig:reg} for the brevity of presentation.
Our results demonstrate that it is essential to select a slightly larger observational
time-step $\dt$ to reduce errors in the estimation of the drift coefficient.
Lasso regression results for the
drift coefficient improve drastically for a larger time-step $\dt=0.01$ compared with 
estimation results for the smaller time-step $\dt=0.001$
(top part of Figure \ref{fig:reg} and Table \ref{table:reg}).
Our numerical results
clearly indicate that selecting a smaller observational time-step has a negative 
effect on the estimation of the drift term since none of the regression techniques can recover the 
correct function form of the drift. This is consistent with our analytical results in 
\eqref{MSEA1} where $\dt$ appears in the denominator.
Polynomial fitting of the diffusion coefficient is comparable for
$\dt=0.001$ and $\dt=0.01$.
Thus, we can see that in contrast with the estimation of the drift term, estimation of the diffusion coefficient is not affected drastically by the observational time-step.
\begin{table}[ht]
\centering
\begin{tabular}{|c|c|c|c|}
\hline
$(M,\dt)$ &  Drift coef.& Diffusion coef. \\
\hline\hline 
True & $-x^3$ & $0.5x^2 +x + 0.5$ \\
\hline
$(1000,0.001)$ &
$-0.8x^7 - 4.3x^5 - 0.002$& 
$0.48x^2 + x + 0.5$\\
\hline
$(1000,0.01)$ &
$-1.03x^3 -0.011x -0.004 $& 
$0.484x^2+x+0.501$\\
\hline
\end{tabular}
\caption{Lasso Polynomial regression fit results for the estimation of the drift and diffusion coefficients for the cubic process \eqref{eq:cubic} with $M=1000$ and $\Delta x=0.05$.}
\label{table:reg}
\end{table}
%
%

\subsubsection*{Numerical Simulations for the Double-Well Potential}
In addition to the model in \eqref{eq:cubic}, we also performed numerical investigation of the double-well potential model with additive and multiplicative noises. In particular, we considered models 
\begin{equation}
\label{dw1}
dX_t = -\gamma X_t (X_t^2 - b_0) + \sigma dW_t,
\end{equation}
and
\begin{equation}
\label{dw2}
dX_t = -\gamma X_t (X_t^2 - b_0) + (\sigma_1 + \sigma_2 X_t^2) dW_t.
\end{equation}
Here the second model is interpreted in \^{I}to sense.
First model is a standard prototype model for many physical processes,
and the choice of diffusion in the second model is motivated by 
a non-trivial diffusion coefficient with non-zero second derivative. 
Parameters in the equations \eqref{dw1} and \eqref{dw2} are 
$\gamma = 2$, $b_0 = 0.5$, $\sigma = \sigma_1 = \sigma_2 = 0.5$.
%
Corresponding stationary variances are approximately $0.42$ and $0.36$.
We performed estimation of the drift and diffusion for equations
\eqref{dw1} and \eqref{dw2}
on intervals $[-1,1]$ and $[-0.6,0.6]$ respectively.

\begin{figure}[ht!]
    \centerline{
    \includegraphics[scale=0.26]{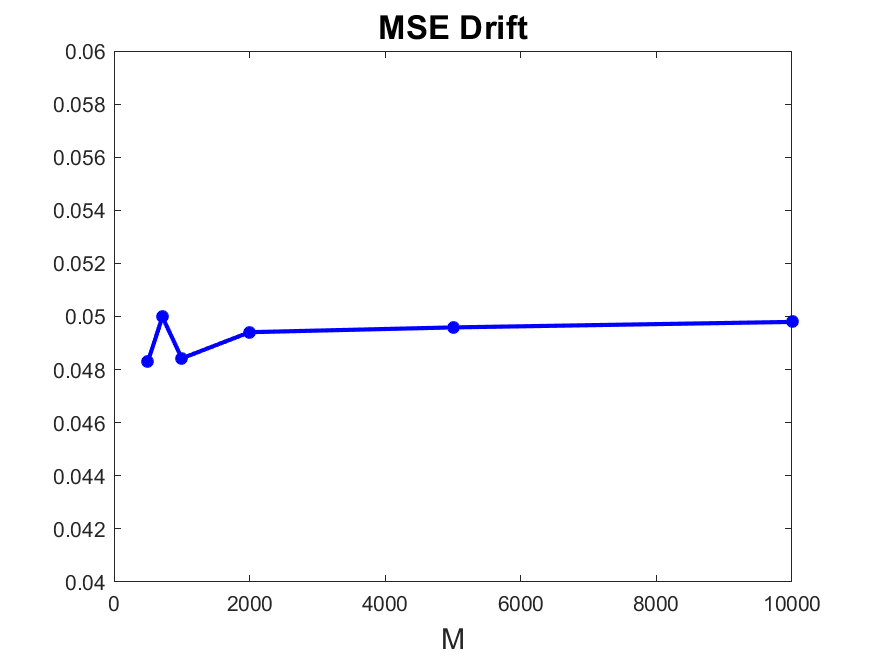} \hspace{-0.4cm}
    \includegraphics[scale=0.26]{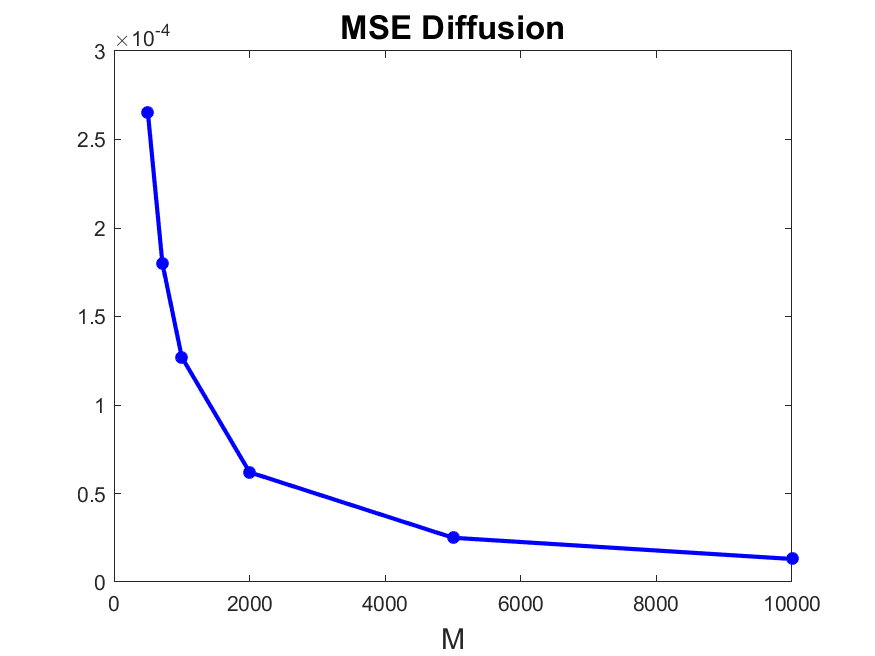}
    \hspace{0.1cm}
    \includegraphics[scale=0.26]{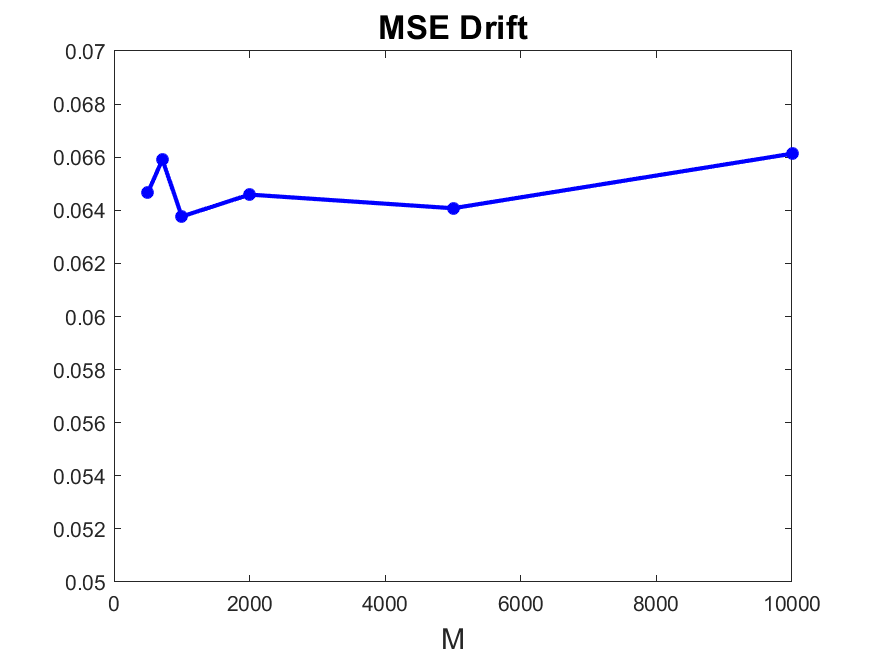} \hspace{-0.4cm}
    \includegraphics[scale=0.26]{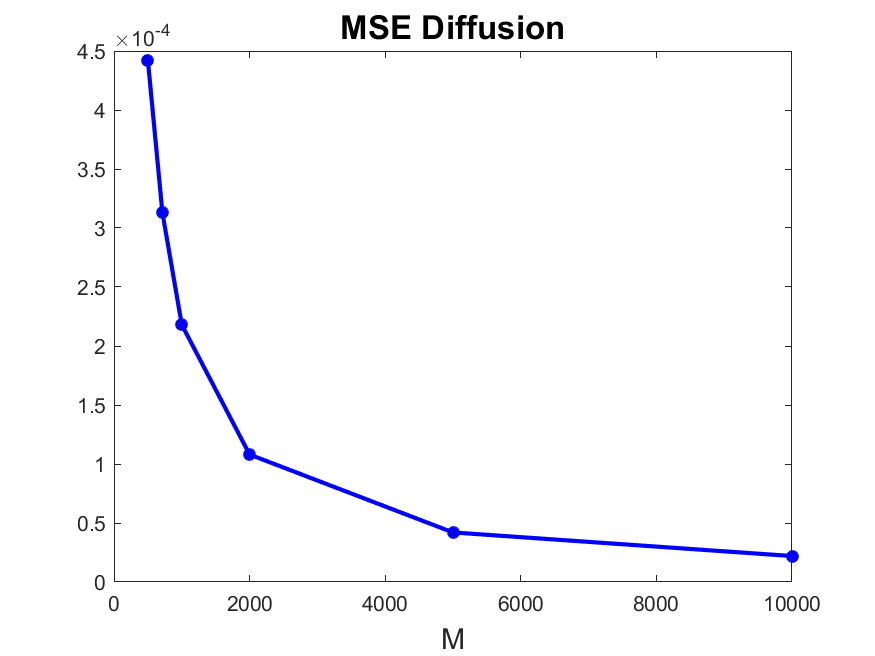}}
    \centerline{
    \includegraphics[scale=0.26]{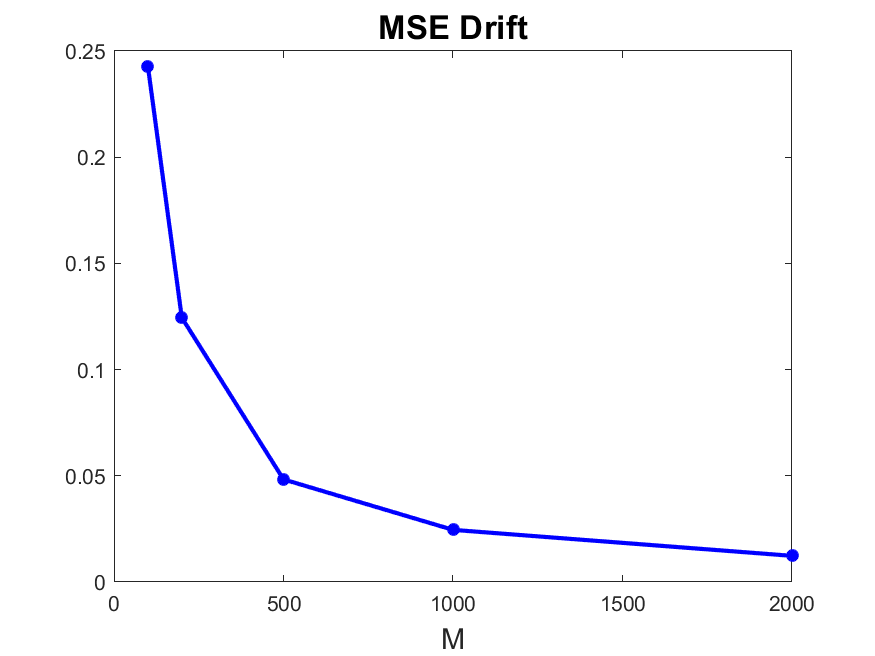} \hspace{-0.4cm}
    \includegraphics[scale=0.26]{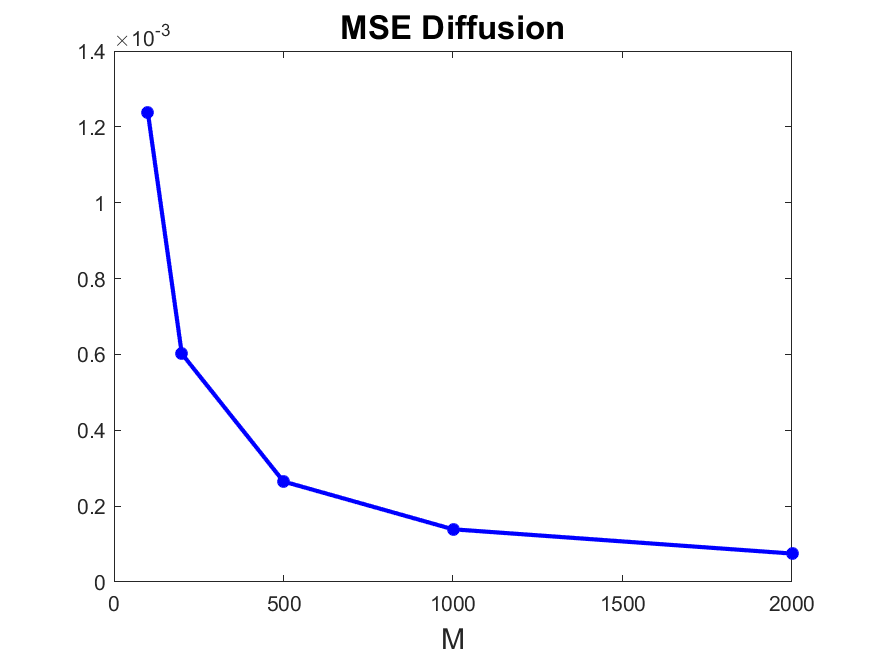}
    \hspace{0.2cm}
    \includegraphics[scale=0.26]{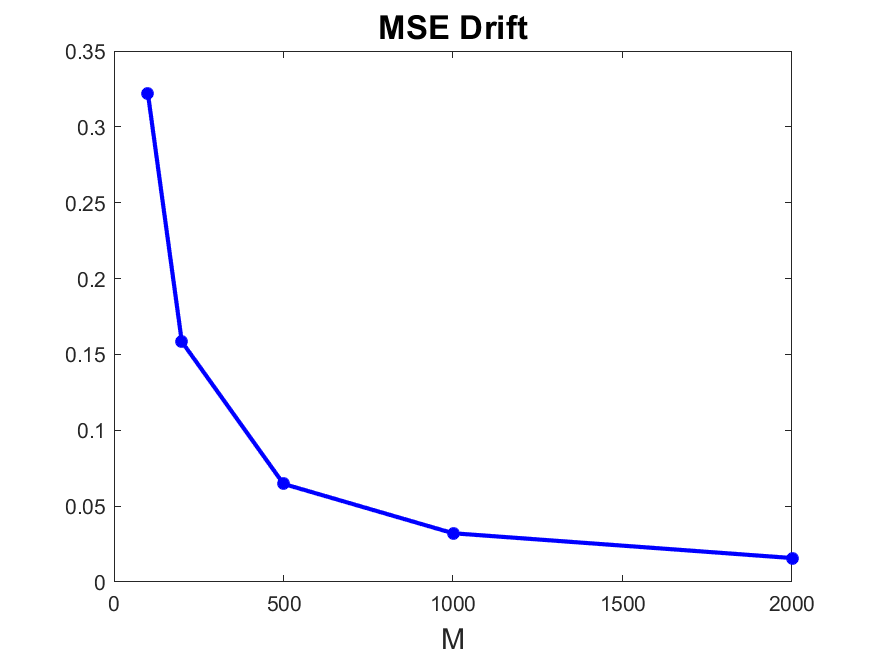} \hspace{-0.4cm}
    \includegraphics[scale=0.26]{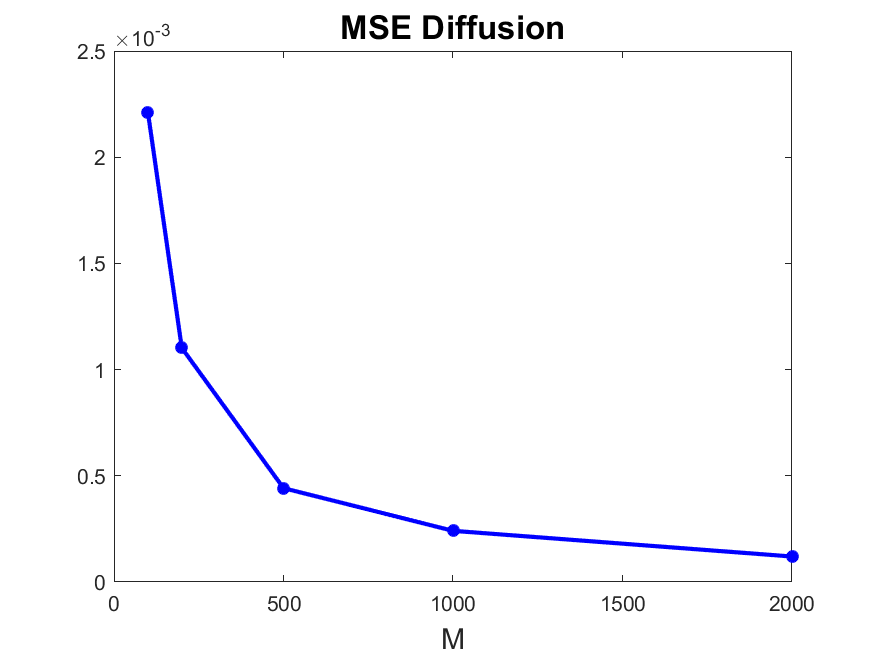}}
    \caption{MSE of the Drift and Diffusion estimators vs $M$
    for equations \eqref{dw1} (two left columns) and \eqref{dw2} (two right columns)
    computed with $\dx=0.1$.
    and two different sampling regimes 
    $M\dt = Const$ (top row) and $M\dt \rightarrow \infty$ (bottom row).
    Corresponding parameters are given by \eqref{mdtconst_dw} and 
    \eqref{mdtinf_dw}.
   }
\label{fig:MSEdw}
\end{figure}

Mean-Squared Error for two different estimation regimes
$M\dt = Cosnt$ and $M\dt \to\infty$ is depicted in Figure \ref{fig:MSEdw}.
Behavior of the Mean-Squared Errors in the estimation of the 
drift and diffusion coefficients for equations \eqref{dw1}
and \eqref{dw2} with respect to two regimes $M\Delta t = Const$ and $M\Delta t\to\infty$ is identical to the estimation of the drift and diffusion for the cubic
equation \eqref{eq:cubic} (c.f. Figures \ref{fig:MSEdw} and \ref{fig:MSE}).
In these simulations we use the following parameters for the two regimes 
in Figure \ref{fig:MSEdw}
\begin{eqnarray}
M\dt &=& 5 \text{~with~} M=500, \, 714, \, 1000, \, 2000, \, 5000, \, 10000,
\label{mdtconst_dw} \\
\dt &=& 0.01 \text{~and~} M= 100, \, 200, \, 500, \, 1000, \, 2000.
\label{mdtinf_dw}
\end{eqnarray}
Similar to the estimation of the cubic process \eqref{eq:cubic}, it is crucial to have large $M\dt$ to reduce errors for the drift estimator. This supports our previous conclusion that it is beneficial 
to perform sub-sampling of the data with a relatively large $\dt$. Estimation of the diffusion 
coefficient improves considerably as $M$ increases and is almost independent of the 
parameter $\dt$. Inspecting our numerical results more closely, we can see that 
estimation errors for the diffusion coefficient are reduced slightly with a smaller $\dt$
(c.f. estimation of the diffusion coefficient in Table \ref{table:dw} with $(M,\dt)=(1000,0.005)$ vs 
$(M,\dt)=(1000,0.01)$ and with $(M,\dt)=(2000,0.0025)$ vs 
$(M,\dt)=(2000,0.01)$).
Therefore, there is a weak dependence on $\dt$, as indicated by the leading order errors terms in \eqref{MSED1}. However, numerical simulations indicate that 
the constant in front of $\dt$ term is smaller compared to the constant in front of 
$M^{-1}$ term. In addition, we we would like to point out that errors for the diffusion 
coefficient are much smaller compared to the errors for the drift.
\begin{table}[ht]
\centering
\begin{tabular}{|c|c|c||c|c|}
\hline
  &
\multicolumn{2}{c||}{Double-Well Model \eqref{dw1}} &
\multicolumn{2}{c|}{Double-Well Model \eqref{dw2}} \\
\hline
$(M,\dt)$ &  Drift coef. & Diffusion coef. & Drift coef. & Diffusion coef. \\
\hline\hline 
$(1000,0.005)$ &  $0.048$ &  $0.00013$    & $0.064$ & $0.00022$\\
\hline
$(1000,0.01)$ &  $0.025$ &  $0.00014$     & $0.032$ & $0.00024$\\
\hline
$(2000,0.0025)$ &  $0.05$ &  $0.000062$    & $0.065$ & $0.00011$\\
\hline
$(2000,0.01)$ &  $0.0123$ &  $0.000075$     & $0.016$  &  $0.00012$\\
\hline
\end{tabular}
\caption{Mean Squared Errors for the estimation of the drift and diffusion coefficients for the double-well models \eqref{dw1} and \eqref{dw2} with $\dx=0.1$, $MC=500$ and different combinations of $(M,\dt)$.}
\label{table:dw}
\end{table}

\begin{figure}[ht!]
\centerline{\includegraphics[scale=0.55]{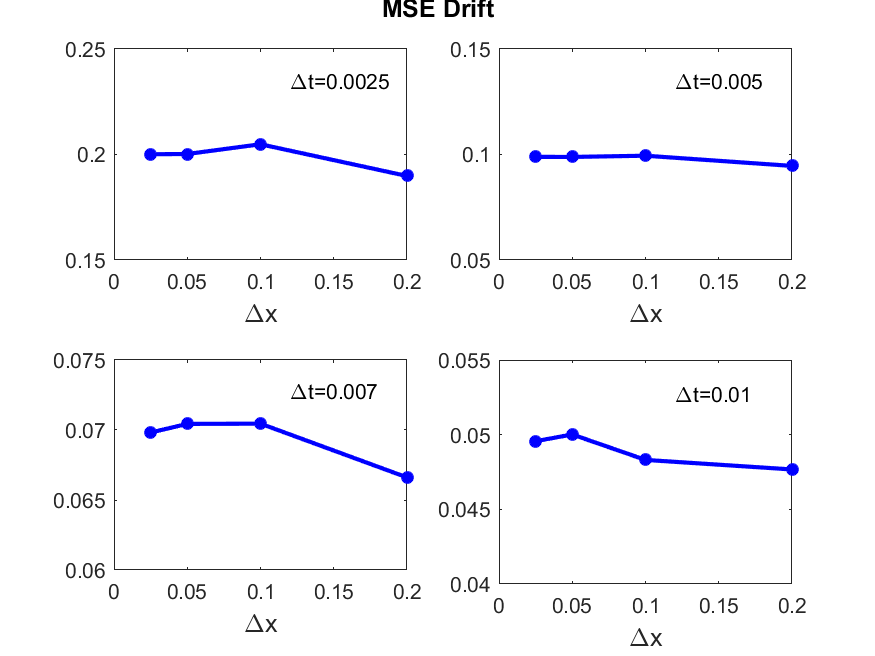}
\hspace*{-0.5cm}
\includegraphics[scale=0.55]{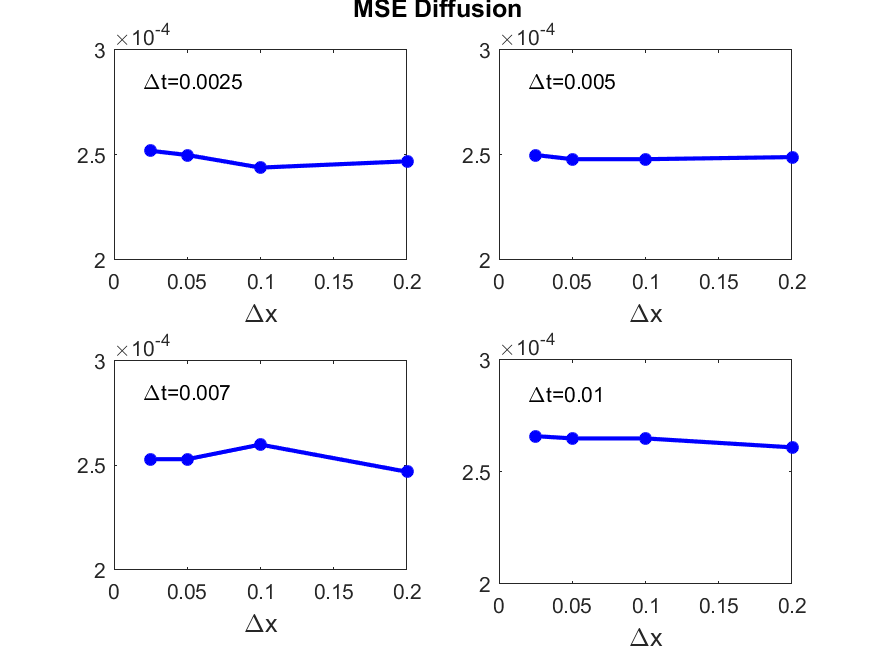}}
\caption{Averaged (over all bins) Mean-Squared-Errors for the estimation of the drift (left) and diffusion (right) coefficients vs $\Delta x$ for several particular values of $\Delta t=0.0025$, $0.005$, $0.007$, $0.01$. Simulations of the double-well process \eqref{dw1} with $M=500$ and $MC=500$.}
\label{fig7}
\end{figure}

\begin{figure}[ht!]
\centerline{\includegraphics[scale=0.55]{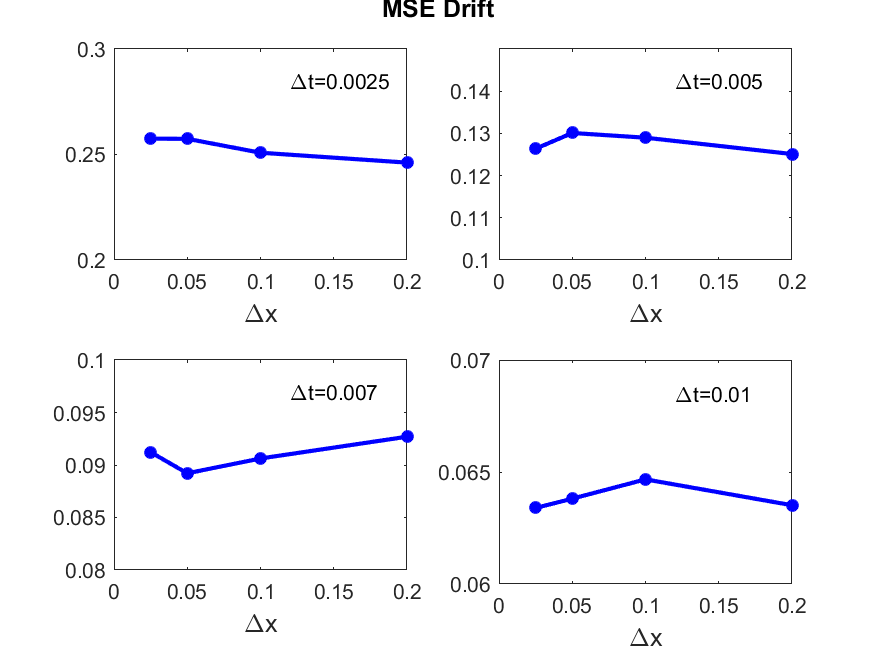}
\hspace*{-0.5cm}
\includegraphics[scale=0.55]{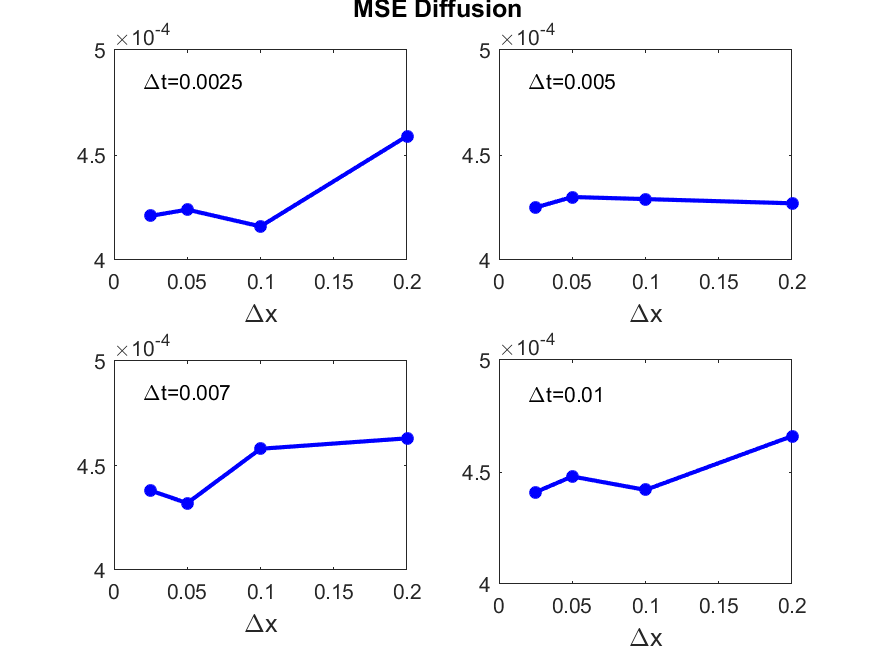}}
\caption{Averaged (over all bins) Mean-Squared-Errors for the estimation of the drift (left) and diffusion (right) coefficients vs $\Delta x$ for several particular values of $\Delta t=0.0025$, $0.005$, $0.007$, $0.01$. Simulations of the double-well process \eqref{dw2} with $M=500$ and $MC=500$.}
\label{fig8}
\end{figure}

Next, Figures \ref{fig7} and \ref{fig8} present dependence of the Mean-Squared Error
in the estimation of the drift and diffusion coefficients for the double-well potential models
\eqref{dw1} and \eqref{dw2} on $\dx$. Similar to the results for the cubic process,
estimators for the drift and diffusion coefficients do not seem to depends strongly on $\dx$
(c.f. Figures \ref{fig7}, \ref{fig8} and Figure \ref{fig:l2snap}).  In addition,
behavior of estimation errors in the $\dx-\dt$ plane for double-well models 
\eqref{dw1} and \eqref{dw2} looks almost identical to Figure \ref{fig:l23d}, except the landscape for the MSE for diffusion estimator is ``more ragged''. 
Thus, we do not to present these results here for the brevity of presentation.
Here we do not present the regression results for the double-well potential models, for the
brevity of presentation. However, it is easy to see that regression results are more accurate 
for the drift estimated with a larger $\dt$, since the men-squared error explodes as $\dt\to 0$ and $M=const$. In addition, estimation of the diffusion coefficient is affected only weakly by using larger values of $\dt$. Finally, we observe almost no dependence of estimation errors on $\dx$
in the parameter regime considered here. 
Overall, numerical results for the double-well models in \eqref{dw1} and \eqref{dw2}
confirm our findings for the cubic model.

\subsubsection*{Discussion}
Our numerical simulations indicate that - 
(i) estimation of both drift and diffusion is not sensitive to the choice of $\dx$
in the parameter regime considered here, 
(ii) estimation of the drift term is affected considerably by the choice of $M\dt$, 
(iii) estimation of the diffusion if primarily affected by $M$ and only slightly by $\dt$.
Therefore, our numerical simulations suggest that for the parameter regime considered here 
the leading sources of error for the drift and diffusion estimation are
\[
\| \hat{A}(x_k) - A(x_k) \|_2^2 \sim \frac{1}{M \dt}, \qquad
\| \hat{D}^2(x_k) - D^2(x_k) \|_2^2 \sim \frac{1}{M}.
\]
Overall, numerical results in this section support our analytical estimates in 
\eqref{MSEA1} and \eqref{MSED1}. However, it is difficult to develop precise 
estimates for the values of constants in front of different terms in \eqref{MSEA1} and \eqref{MSED1}. It is possible to obtain a different asymptotic behavior by considering datasets with very large values of $M$. However, our numerical simulations indicate that practical values $M=500,\ldots,1000$ 
are sufficient for accurate estimation of drift and diffusion coefficients in many situations.

Practical applications of the estimation approach discussed here include two distinct cases - (i) when dataset for estimation is generated by a numerical simulation of a complex model (e.g. turbulence) and (ii) when dataset is generated by a physical experiment.
In the first case we have a full control over the estimation parameters $M$, $\dx$, and $\dt$, but in the second case the choice of 
estimation parameters can be limited by the experimental setup. Numerical simulations presented here indicate that it is sufficient to have
approximately 10 - 20 bins over the interval $[\text{mean} \pm \text{standard deviation}]$. This number of bins is sufficient for applying regression techniques to the estimated drift and diffusion coefficients if they are relatively slow varying (e.g. polynomials of order 3 - 5).
One can start by applying the estimation procedure with 10 bins and verify the result with 15 and 20 bins. However, we would like to 
point out that this number of bins might be insufficient for estimating the drift and diffusion coefficients from multiscale data 
(e.g. highly oscillatory potential). Thus, dealing with multiscale problems requires generating considerably larger datasets. 
From our experience with numerical simulations, it is sufficient to have the number of data points in each bin of the order $M=500,\ldots,1000$ and
even estimation with as low as $M=200$ points can produce adequate results. 

Selecting appropriate  time-step of estimation $\dt$ 
is probably the central issue for obtaining good estimation results.
Optimal choice of $\dt$ is most likely related to the time-scale of the dynamic variables. 
In the numerical example presented in this paper, the correlation time of the process \eqref{eq:cubic} with the chosen parameters is 
$CT \approx 2$. Therefore, choice of $\dt \approx [0.005,\ldots,0.02]$ seems appropriate for data with correlation time of $O(1)$.
Computational challenges can arise when estimating multi-dimensional multiscale systems where dynamic variables have considerably different
correlation times. 
A practical guideline for verifying numerical estimation results should be computing estimates 
with the time-step $\dt$ and with a larger time-step $2\dt$ or even $4\dt$. Then, estimation with a larger $\dt$ can be taken as the'' truth''
in the computations of the MSE errors with smaller time-step. Small practical errors for the drift estimator are of order $O(10^{-1})$.

Many practical situations require post-processing of previously generated datasets.
In such situations, the value of $M$ (number of points in each bin) is not fixed, but varies from one bin to another. Of course, in such cases all available points should be utilized for estimation in a particular bin. Moreover, our numerical tests indicate that since errors do no depend on the bin size for a wide range of computational parameters, bins should be taken quite large to ensure
better estimation results.
However, a situation might arise when different bins have very different number
of estimation points. For instance, bins in the tails of the stationary distribution or near the top of a well (i.e. near $x=0$) in a double-well potential might have very few estimation points.
In such cases, one should neglect bins with a small number of 
estimation points when applying nonlinear regression. A practical guideline is to neglect bins with fewer than $M=200$ or $M=500$ points, depending on the availability of data in other bins.
Thus, one can set a threshold for the number of estimation points in each bin to ensure 
better result for nonlinear regression, since nonlinear regression can be easily applied on a non-uniform mesh.

\section{Conclusion}
\label{sec:conc}
In this paper we establish rigorous mathematical foundation for the optimal choice of computational 
parameters for estimators of the drift and diffusion coefficients in stochastic differential equations 
driven by Brownian motion based on conditional expectations. 
It has become viable to utilize this approach for higher-dimensional problems due to increase in computational capabilities and 
increasing availability of data. However, an important practical task is 
to optimize selection of computational parameters in order to minimize computational and data-generating complexities.
To address this issue,  
we analyze asymptotic behavior of the bias and mean squared error for both estimators and arrive at important practical results for the selection of computational sampling parameters. In particular, 
we demonstrate both analytically and numerically that the spatial mesh size for estimation can be taken
much larger than the observational time-step and the approximate practical scaling 
of space- and time-discretization parameters should be
\[
\dx \sim \sqrt{\dt}.
\]
This scaling has important practical implications, especially for higher-dimensional problems. 
In particular, this scaling implies that in many practical applications
the space-discretization (bin size) should be taken much larger than the time-discretization. This significantly 
reduces the computational complexity of the problem since the bin size can be taken to be quite large. Large bin sizes also
imply that the observational trajectory used for estimation can be short, observational points are more likely to ``fill-up''
bins of a larger size. In addition, one can potentially develop estimation strategies with \emph{overlapping bins} where
the one observation point would contribute to the estimation for two neighboring spatial points $x_k$ and $x_{k+1}$
such that $x_{k+1} - x_k < \dx$.  This approach can be used to reduce estimation errors due to a small sample size, $M$.

Overall, our numerical simulations support analytical expressions in 
 \eqref{MSEA1} and \eqref{MSED1}. However, since numerical errors are not sensitive to the changes in the bin size, $\dx$, 
Our numerical simulations indicate that the leading terms in the drift and diffusion estimation are
\begin{equation}
\label{MSEAD}
\| \hat{A}(x_k) - A(x_k) \|_2^2 \sim \frac{1}{M \dt}, \qquad
\| \hat{D}^2(x_k) - D^2(x_k) \|_2^2 \sim \frac{1}{M}.
\end{equation}
Absence of influence of other terms in the numerical errors for the drift and diffusion  
(c.f. \eqref{MSEA1} and \eqref{MSED1} with expressions above) can be attributed to small 
constants in front on other asymptotic terms in 
 \eqref{MSEA1} and \eqref{MSED1}.
 
Spatial refinement with respect to the bin size does not yield significant improvement for both, the drift and the
diffusion estimators. This suggests that estimation errors due to the space-discretization are typically much smaller, compared to other sources of error. However, balance between errors due to temporal and spacial
discretizations might depend on the roughness of the data. For instance, for SDEs with highly-oscillatory multi-scale 
potential errors due to the spatial discretization might play an important role. 
This will be examined in a subsequent paper.

Analytical and numerical results presented in this paper provide guidelines for developing practical estimation
schemes for the drift and diffusion coefficients from stationary time-series. In particular, estimation with 
$M=500,\ldots,1000$, $\dx=0.05,\ldots,0.1$, $\dt \approx 0.01$ provides good results for both, the drift and the diffusion.
Although here we only 
consider a scalar SDE, we expect that conclusions reached in this paper should hold for systems of equations as well.
We will verify this numerically in a subsequent paper.

General definitions of the drift and diffusion in \eqref{Ax} and \eqref{Bx} are applicable to non-stationary time-series as well. However, the practical application of these formulas for stationary and non-stationary time-series differs considerably. In particular, for stationary time-series we use ergodicity of the stochastic process $X_t$. For non-stationary processes, a single individual trajectory can drift to infinity and, thus, there might be very few points in any particular bin. This implies that in practice estimation in a non-stationary regime should utilize ensemble simulations or repeated experiments 
and index sets $M_k$ in \eqref{Ahat}, \eqref{Dhat} should use data points from multiple trajectories. 
For non-stationary time-series, general expressions for the asymptotic behavior of errors in  
\eqref{MSEA1} and \eqref{MSED1} remains the same. In practice, there might additional numerical cancellations
of these errors due to uncorrelated data in ensemble trajectories. However, since our numerical simulations indicate that 
the leading behavior of MSEs
for the drift and diffusion is given by \eqref{MSEAD} in many practical situations, we expect the the same behavior when 
estimates \eqref{Ahat}, \eqref{Dhat} are applied to non-stationary ensemble data.

Finally, the non-parametric estimation framework discussed here can also be combined with
regression techniques to perform parametric fitting of a nonlinear function to obtain the functional form of the
drift and diffusion coefficients (e.g. \cite{clementi2}). In addition, it is also possible to utilize LASSO-type techniques
to obtain the optimal functional form of the drift and diffusion. Future research directions will focus on these
issues with the emphasis on the practical aspect of estimating drift and diffusion coefficients for systems of
SDEs.

\appendix
\section{MSE of the Drift Estimator}
\label{ap1}

Here we outline the calculation for the MSE of the drift estimator $\hat{A}(x_k)$.
As discussed in section \ref{sec:msedrift}
\begin{align}
& MSE\{\hat{A}(x_k)\}  \approx
\nonumber \\
& \bE \left[ \left(\frac{1}{M\dt}   \sum\limits_{j \in M_k} \left( [A(X_{t_j}) - A(x_k)]\dt + \sum\limits_{q=1}^6 B_q(X_{t_j})I_{\alpha_q, j} \right)\right)^2 \Big| X_{t_j}\in Bin_k \right] =
\nonumber \\
& \bE \left[ \frac{1}{M^2\dt^2}   \sum\limits_{i,j \in M_k}\left(  [A(X_{t_j}) - A(x_k)]\dt + \sum\limits_{q=1}^6 B_q(X_{t_j})I_{\alpha_q, j} \right) \right. \times 
\nonumber \\
& \qquad \qquad \qquad \qquad 
\left(  [A(X_{t_i}) - A(x_k)]\dt + \sum\limits_{l=1}^6 B_l(X_{t_i})I_{\alpha_l, i} \right)
\Big| X_{t_i}, X_{t_j}\in Bin_k \Bigg].
\label{expr1}
\end{align}
Here, we proceed by expanding the square, but we keep the term $[A(X_{t_j}) - A(x_k)]$ together. We need to consider double summations with
products of stochastic integrals with different indexes.
We outline below types of terms which are treated differently and point out leading 
order terms for each type.

\noindent
{\bf{Type 1}:} Consider the cross-product of the first two terms in \eqref{expr1}
\begin{align*}
&\frac{1}{(M\dt)^2 }\sum_{i, j \in M_k}  \labs \bE \left[ (A(X_{t_i}) - A(x_k)) \, (A(X_{t_j}) - A(x_k)) \dt^2
\labs X_{t_i}, X_{t_j} \in Bin_k \right. \right]  \rabs \\
&\quad \le \frac{1}{M^2 }\sum_{i, j \in M_k}  \bE \left[ K_A^2 \labs X_{t_i} - x_k \rabs \labs X_{t_j} - x_k \rabs
\Bigg|  X_{t_i}, X_{t_j} \in Bin_k  \right] \\
& \quad \le C \dx^2,
\end{align*}
where we used that $A(x)$ is Lipschitz and $|X_{t_j} - x_k|, \, |X_{t_i} - x_k| \le \dx/2$ since both
$X_{t_j}, \, X_{t_i} \in Bin_k$ and $x_k$ is the center of the bin.

\noindent
{\bf{Type 2}:} Consider cross-terms of the form 
\begin{eqnarray*}
&& \frac{1}{(M\dt)^2} \sum_{i, j \in M_k} \bE \left[
(A(X_{t_i}) - A(x_k)) \dt \, B_q(X_{t_j})I_{\alpha_q, j} + \right. \\
&& \qquad \qquad \qquad \left. (A(X_{t_j}) - A(x_k)) \dt \, B_l(X_{t_i})I_{\alpha_l, i}
\Bigg| X_{t_i}, X_{t_j} \in Bin_k  \right]  \\
&& = \frac{2}{M^2 \dt} \sum_{i, j \in M_k} \bE \left[(A(X_{t_i}) - A(x_k)) \, B_q(X_{t_j})I_{\alpha_q, j} 
\Bigg| X_{t_i}, X_{t_j} \in Bin_k  \right] ,  \text{~~for~~} q=1,\ldots,6
\end{eqnarray*}
where we used symmetry between $t_i$ and $t_j$ and $\alpha_l$ and $\alpha_q$.
Here $X_{t_i}$ and $I_{\alpha_l, j}$
are \emph{not} independent if $t_i > t_j$.
Therefore, we use the Lipschitz property of $A(x)$ and obtain
\begin{align*}
& \frac{2}{M^2 \Delta t}  \labs \bE\left[ \sum\limits_q  \sum_{i, j \in M_k}(A(X_{t_i})-A(x_k))
B_q(X_{t_j})I_{\alpha_q, j} \Bigg| X_{t_i}, X_{t_j} \in Bin_k  \right]  \rabs \\
& \quad \le \frac{2}{M^2\dt} \sum\limits_q \sum_{i, j \in M_k}\bE \left[ \labs \left(A(X_{t_i}) - A(x_k) \right)
B_q(X_{t_j})I_{\alpha_q, j} \rabs \Bigg| X_{t_i}, X_{t_j} \in Bin_k  \right]  \\
& \quad \le \frac{2K_A\dx}{M^2 \dt} \sum\limits_q \sum_{i, j \in M_k}\bE \left[ \labs 
B_q(X_{t_j})I_{\alpha_q, j} \rabs \Bigg| X_{t_j} \in Bin_k  \right]  \\
& \quad \le \frac{2K_A \Delta x}{M^2 \Delta t} \sum\limits_q \sum_{i, j \in M_k}
\left( \bE_k B_q^2(X_{t_j}) \right)^{1/2} 
\left( \bE I_{\alpha_q, j}^2 \right)^{1/2} \\
& \quad \le \frac{C\dx}{\sqrt{\dt}} \left( 1 + \sqrt{\dt} + O(\dt^{3/2}) \right),
\end{align*}
where we used the H\"{o}lder inequality 
and lowest-order terms are due to 
$\bE I_{(1), j}^2 =\dt $. 
Other stochastic integrals contribute to higher-order terms.
Here we use a notation for the conditional expectation 
$\bE_k f(x) = \bE [f(x) | x \in Bin_k]$.
Since the truncated density has a finite support, we assume that all conditional expectations exist and are finite., e.g.,
$\bE_k B_q^2(X_{t_j}) < \infty$.

\noindent
{\bf{Type 3}:} Consider terms with stochastic integrals for either $q = 5$ or $l = 5$
\begin{align*}
& \alpha_q = (0,0) \text{~~and~~} \alpha_l = (1), (1,1), (1,0), (0,1), (0,0), (1,1,1), \\
& \alpha_l = (0,0) \text{~~and~~}  \alpha_q = (1), (1,1), (1,0), (0,1), (1,1,1).
\end{align*}
Due to symmetry, we only need to consider $q = 5$. Recall that $\alpha_5=(0,0)$ and 
$I_{(0,0),j} = \dt^2/2$. 
Then \eqref{expr1} becomes
\begin{eqnarray*}
&& \frac{1}{M^2} \labs
\sum\limits_{i, j \in M_k} \bE \left[B_5(X_{t_i})  B_q(X_{t_j}) I_{\alpha_q, j} \Bigg| X_{t_i}, X_{t_j} \in Bin_k  \right]  \rabs \le
 \\
&&  \frac{1}{M^2} \sum\limits_{i, j \in M_k} 
\left( \bE \left[B_5^2(X_{t_i}) \,  B_q^2(X_{t_j})  \left| X_{t_i}, X_{t_j} \in Bin_k\right. \right] \right)^{1/2} \| I_{\alpha_q, j} \|_2,
\end{eqnarray*}
where used the H\"{o}lder inequality. 
Similar to Type 2 terms, we assume that all expectations with respect to the joint 
truncated density exist and are finite
(i.e. $\bE_k \left[B_2^2(X_{t_i}) \,  B_q^2(X_{t_j}) \right] < \infty$). The exact form of this joint density is hard to analyze, 
but it has a finite support and, thus, this assumption is quite reasonable.

As a final step, we only need to analyze lowest-order terms resulting from stochastic integrals. 
Therefore, 
\[
\frac{1}{2M^2}  \labs \sum_q \sum_{i, j \in M_k}  \bE \left[B_5(X_{t_i}) \,  B_q(X_{t_j})I_{\alpha_q, j} \labs X_{t_i}, X_{t_j} \in Bin_k\right.\right] \rabs \le 
 C \sqrt{\dt} \left(1 + \sqrt{\dt} + O(\dt^{3/2}) \right),
\]
where the lowest-order term is due to $\| I_{(1), j} \|_2 = \sqrt{\dt}$ and the next term arises from
$\| I_{(1,1), j} \|_2 = \dt/\sqrt{2}$. Here $C$ is some generic constant representing 
upper bound for all expectations of the form $\bE_k \left[B_5^2(X_{t_i}) \,  B_l^2(X_{t_j})\right]$.

\noindent
{\bf{Type 4}:} Consider all possible combinations of stochastic integrals 
with the following indexes
\begin{equation}
\alpha_q, \alpha_l = (1), (1,1), (0,1), (1,0), (1,1,1).\label{msedrift2}
\end{equation}
Without loss of generality, we can assume that $t_j > t_i$. 
Then from the property of stochastic integrals
\[
\bE \left[B_q(X_{t_i})I_{\alpha_q, i} B_l(X_{t_j})I_{\alpha_l, j} \right] = 
\bE \left[B_q(X_{t_i})I_{\alpha_q, i} B_l(X_{t_j})\right] \bE \left[ I_{\alpha_l, j} \right] = 0.
\]
Therefore, for these combinations of stochastic integrals we only need to 
consider case $i = j$ and terms \eqref{expr1}
become
\begin{align*}
& \frac{1}{(M\dt)^2} \labs\sum_{i, j \in M_k} \bE \left[B_q(X_{t_i}) I_{\alpha_q, i} \, B_l(X_{t_j}) I_{\alpha_l, j} \labs  X_{t_i}, X_{t_j} \in Bin_k \right. \right]\rabs   = \\
& \frac{1}{(M\dt)^2} \labs \sum_{i \in M_k} \bE \left[B_q(X_{t_i}) \, B_l(X_{t_i}) I_{\alpha_q, i} I_{\alpha_l, i} \labs X_{t_i} \in Bin_k \right. \right]\rabs \le \\
& \frac{1}{(M\dt)^2} \sum_{i \in M_k} \left( \bE_k \left[ B_q^2(X_{t_i}) \, B_l^2(X_{t_i}) 
\right]\right)^{1/2} \| I_{\alpha_q, i} I_{\alpha_l, i} \|_2.
\end{align*}
Here we need to calculate 
fourth moments of stochastic integrals. The lowest-order term is due to 
$\| I_{(1), i}^2 \|_2 = \sqrt{3} \dt $. All other combinations of stochastic integrals result in 
moments of higher order (some of them are given by 5.2 and 5.7 of \cite{Peter}), e.g.
\begin{eqnarray*}
 && \| I_{(1), i} I_{(1, 1), i} \|_2 = O(\dt^{3/2}), \qquad 
 \| I_{(1, 1), i}^2 \|_2 = \frac{\sqrt{15}}{2}(\dt)^2,\\
 && \| I_{(0, 1), i}^2 \|_2 = \| I_{(1, 0), i}^2 \|_2 = \| I_{(1, 1,1), i}^2 \|_2 = O(\dt^3).
\end{eqnarray*}
Therefore, 
\[
\frac{1}{(M\dt)^2} \labs \sum\limits_{q,l} \sum_{\substack{i \in M_k}} \bE \left[B_q(X_{t_i}) \, B_l(X_{t_i}) I_{\alpha_q, i} I_{\alpha_l, i} \labs X_{t_i} \in Bin_k \right. \right] \rabs
\le \frac{C}{M \dt} \left( 1 + \sqrt{\dt} + O(\dt) \right), 
\]
where summation with respect to $q,l$ is taken over \eqref{msedrift2} and 
$C$ is a suitable constant.

\section{MSE of the Diffusion Estimator}
\label{ap2}

In this section, we focus our attention on the MSE of diffusion estimator given by \eqref{MSEDDEF} in section \ref{sec:msediff}.
The MSE squared of the diffusion estimator is given by
\small
\begin{align}
& \|\hat{D}^2(x_k) - D^2(x_k) \|_2^2 \approx
\bE \left[
\left( \frac{1}{M\dt}    \sum\limits_{j \in M_k} 
\sum\limits_{l, q=0}^6 B_l(X_{t_j}) B_q(X_{t_j}) I_{\alpha_l, j} I_{\alpha_q, j}  - D^2(x_k) \right)^2
\Big| X_{t_j} \in Bin_k \right] = 
\nonumber \\
%
& 
\underbrace{ 
\frac{1}{(M\dt)^2}  \bE \left[
\sum\limits_{i,j \in M_k} 
\left( D^2(X_{t_i}) I_{(1), i}^2 - D^2(x_k) \dt \right)
\left( D^2(X_{t_j}) I_{(1), j}^2 - D^2(x_k) \dt \right) \Big| X_{t_i}, X_{t_j} \in Bin_k \right] }_{Type 1} + 
\label{term1} \\ 
%
& 
\underbrace{ 
\frac{1}{(M\dt)^2}  \bE \left[ \left( \sum\limits_{j \in M_k} \sum\limits_{\substack{l, q=0 \\ l \times q \ne 1}}^6 B_l(X_{t_j}) B_q(X_{t_j}) I_{\alpha_l, j} I_{\alpha_q, j} \right)^2
\Big| X_{t_j} \in Bin_k  \right] }_{Type 2} +
\label{term2} \\
%
& \underbrace{
\frac{2}{(M\dt)^2} \bE \left[  \left( \sum\limits_{i,j \in M_k} 
\left( D^2(X_{t_i}) I_{(1), i}^2 - D^2(x_k) \dt \right)
\sum\limits_{\substack{l, q=0 \\ l \times q \ne 1}}^6 B_l(X_{t_j}) B_q(X_{t_j}) I_{\alpha_l, j} I_{\alpha_q, j} \right)
\Big| X_{t_i}, X_{t_j} \in Bin_k \right] }_{Type 3} .
\label{term3}
\end{align}
\normalsize
We proceed by considering the terms above separately.

\noindent
{\bf{Type 1}:} 
First, we consider the first term in \eqref{term1} and
by adding and subtracting $D^2(X_{t_i})\dt$ and $D^2(X_{t_j})\dt$ in the 
first and second bracket, respectively, we obtain
\begin{eqnarray*}
&&\frac{1}{(M\dt)^2}  \bE \left[
\sum\limits_{i,j \in M_k}
\left( D^2(X_{t_i}) I_{(1), i}^2 - D^2(x_k) \dt \right)
\left( D^2(X_{t_j}) I_{(1), j}^2 - D^2(x_k) \dt \right)\right] = \\
&& \quad \frac{1}{(M\dt)^2}  \bE \left[
\sum\limits_{i,j \in M_k} D^2(X_{t_i})D^2(X_{t_j}) \left(I_{(1), i}^2 -\dt \right)\left(I_{(1), j}^2 -\dt \right)\right] + \\
&& \quad  \frac{2}{M^2\dt}  \bE \left[ \sum\limits_{i,j \in M_k}  D^2(X_{t_i}) \left(I_{(1), i}^2 -\dt \right)  \left( D^2(X_{t_j}) - D^2(x_k) \right)\right] + \\
&& \quad  \frac{1}{M^2}  \bE \left[ \sum\limits_{i,j \in M_k}  \left( D^2(X_{t_i}) - D^2(x_k) \right) \left( D^2(X_{t_j}) - D^2(x_k) \right) \right] \le \\
&& \frac{1}{(M\dt)^2} \sum\limits_{i \in M_k} \bE_k \left[D^4(X_{t_i}) \right]\bE \left[\left(I_{(1), i}^2 -\dt \right)^2 \right] +
 C\frac{K_D \dx}{\dt} \bE \left[ |I_{(1), i}^2 -\dt |\right] + 
 \frac{(K_D \dx)^2}{4}  = \\
&& C \left(\frac{1}{M} +  \dx + \dx^2\right),
\end{eqnarray*}
where we used that 
$\bE \left[\left(I_{(1), i}^2 -\dt \right)^2 \right] = O(\dt^2)$,
$\bE \left[ \left| I_{(1), i}^2 -\dt \right| \right] = O(\dt)$, 
$K_D$ is a Lipschitz constant for $D^2(x)$,
and we use $C$ to denote some generic constant.

\noindent
{\bf{Type 2}:} Consider the terms arising from \eqref{term2}. There are a lot of terms arising from squaring the sum in  \eqref{term2}, but all of them have the following form
\begin{equation}
\label{expr2}
\frac{1}{(M\dt)^2}  \bE \left[\sum\limits_{i,j \in M_k} 
B_q(X_{t_i}) B_l(X_{t_i}) I_{\alpha_q, i} I_{\alpha_l, i}
B_r(X_{t_j}) B_m(X_{t_j}) I_{\alpha_r, j} I_{\alpha_m, j}
 \right],
\end{equation}
where $q, l, r, m = 0, \ldots, 6$ with restriction $q\times l \ne 1$ and $r \times m \ne 1$ since
cases $q\times l = 1$ and $r \times m = 1$ correspond to terms of type 1 and type 3 considered separately. This means that for type 2 terms we cannot have $\alpha_q = \alpha_l = (1)$ or
$\alpha_r = \alpha_m = (1)$.
Here indexes $q,l$ correspond to time $t_i$ and indexes $r,m$ correspond to time $t_j$.
There are many terms of type 2 and we will distinguish several sub-types. 

\noindent
{\bf{Type 2a}:} \\
Consider type 2 terms with the following 3 restrictions - \\ 
(i) at least one of the integrals in each pair $I_{\alpha_q, i} I_{\alpha_l, i}$ and $I_{\alpha_r, j} I_{\alpha_m, j}$ is stochastic,\\
(ii) $q \ne l$ and $r \ne m$,\\
(iii) $n(\alpha_q, \alpha_l)$ and $n(\alpha_r, \alpha_m)$ (number of 1's) are odd.

Without loss of generality we can consider $t_j > t_i$.
Then we can write
\begin{eqnarray*}
&& \bE \left[
B_q(X_{t_i}) B_l(X_{t_i}) I_{\alpha_q, i} I_{\alpha_l, i}
B_r(X_{t_j}) B_m(X_{t_j}) I_{\alpha_r, j} I_{\alpha_m, j} \right] =
\\
&& \bE \left[
A_q(X_{t_i}) B_l(X_{t_i}) I_{\alpha_q, i} I_{\alpha_l, i}
A_r(X_{t_j}) B_m(X_{t_j}) \right]
\bE \left[ I_{\alpha_r, j} I_{\alpha_m, j} \right]
\end{eqnarray*}
since time intervals $[t_i, t_i+\dt]$ and $[t_j, t_j+\dt]$ do not overlap.
And using the fact that the number of 1's in the pair $(\alpha_r, \alpha_m)$ is odd, 
\[
\bE \left[ I_{\alpha_r, j} I_{\alpha_m, j} \right] = 0.
\] 
where we use Lemma 5.7.2 in \cite{Peter}. A similar argument holds for $t_i > t_j$.
Therefore,  due to the condition (iii), we can reduce the 
sum in \eqref{expr2} to the case $i=j$, i.e.,
\begin{eqnarray*}
&& \frac{1}{(M\dt)^2}  \bE \left[\sum\limits_{i,j \in M_k} 
B_q(X_{t_i}) B_l(X_{t_i}) I_{\alpha_q, i} I_{\alpha_l, i}
B_r(X_{t_j}) B_m(X_{t_j}) I_{\alpha_r, j} I_{\alpha_m, j}
 \right]  = \\
&& \frac{1}{(M\dt)^2}  \bE \left[\sum\limits_{j \in M_k} 
B_q(X_{t_j}) B_l(X_{t_j}) I_{\alpha_q, j} I_{\alpha_l, j}
B_r(X_{t_j}) B_m(X_{t_j}) I_{\alpha_r, j} I_{\alpha_m, j}
 \right].	
\end{eqnarray*}
This reduces the number of terms in the summation from $M^2$ (when $i,j \in M_k$)
to $M$ (when $j\in M_k$).

There are many combinations of  indexes $\alpha_q, \alpha_l, \alpha_r, \alpha_m$ 
which satisfy requirements for type 2a terms.
Since $I_{(1), j} \sim \sqrt{\dt}$
are lowest-order stochastic integrals in the expansion \eqref{ITE}, lowest-order terms
for type 2(i) will appear when $q = r = 1$ (or when  $l = m = 1$ by symmetry). In this case
$\alpha_q = \alpha_r = (1)$.
Due to the restriction (ii) for type 2a terms, if $q = r = 1$, then $l \ne 1$ and $m \ne 1$. 
In addition, it is also clear that in order to capture the leading-order type 2a terms
integrals $I_{\alpha_l,i}$ and $I_{\alpha_m,j}$ should be of the lowest possible order.
There are two integrals of order $\dt$, namely $I_{(0),i} = \dt$ and $I_{(1,1),i} \sim \dt $
(since $||I_{(1,1),i}||_2 \sim \dt$).
Therefore, to obtain lowest-order type 2a terms, 
indexes $\alpha_l$ and $\alpha_m$ should correspond to those two integrals.

Thus, here we list some lower order terms of Type 2a. \\
(a) $q = r = 1$ and $l = m = 2$ or we can switch $q$, $l$ and $r$, $m$ because of symmetry.
In this case $\alpha_q = \alpha_r = (1)$ and $\alpha_l = \alpha_m = (1,1)$.
\begin{align*}
& \frac{1}{(M\dt)^2}  \sum\limits_{i,j \in M_k}   \bE \left[ 
B_1(X_{t_i}) B_2(X_{t_i}) I_{(1), i} I_{(1,1), i}
B_1(X_{t_j}) B_2(X_{t_j}) I_{(1), j} I_{(1,1), j}
 \right] = \\
& \frac{1}{(M\dt)^2}   \sum\limits_{j \in M_k}  \bE \left[
B_1^2(X_{t_j}) B_2^2(X_{t_j}) I_{(1), j}^2 I_{(1,1), j}^2
\right] = \\
& \frac{1}{M\dt^2} \bE_k\left[B_1^2(x) B_2^2(x)\right] 
\| I_{(1), j}^2 I_{(1,1), j}^2 \|_2^2 
\le \dfrac{C\dt}{M}.
\end{align*}
(b) $q = r = 1$ and $l = m = 0$ or we can switch $q$, $l$ and $r$, $m$.
In this case $\alpha_q = \alpha_r = (1)$ and $\alpha_l = \alpha_m = (0)$ and
we would like to remind that $I_{(0), i} = I_{(0), j} = \dt$. 
Therefore, \eqref{expr2} reduces to 
\begin{align*}
& \frac{1}{M^2}   \sum\limits_{i,j \in M_k}  \bE \left[ 
B_1(X_{t_i}) B_0(X_{t_i}) I_{(1), i} 
B_1(X_{t_j}) B_0(X_{t_j}) I_{(1), j} 
 \right]  = \\
& \frac{1}{M^2}   \sum\limits_{j \in M_k}  \bE \left[
B_1^2(X_{t_j})B_0^2(X_{t_j}) I_{(1), j}^2 
 \right]   = \\
& \frac{1}{M^2}   \sum\limits_{j \in M_k}  \bE_k \left[
B_1^2(X_{t_j})B_0^2(X_{t_j})\right] \bE \left[ I_{(1), j}^2 
 \right] 
  \le \frac{C\dt}{M}.
\end{align*}
(c) $q =r =1$, $l = 0$ and $m = 2$ 
 or we can switch $q$, $l$ and $r$, $m$.
In this case $I_{\alpha_l, i} = I_{(0), i} = \dt$, $I_{\alpha_m,j} = I_{(1,1),j}$ and, therefore, 
\eqref{expr2} becomes
\begin{align*}
& \frac{1}{(M\dt)^2}  \sum\limits_{i,j \in M_k}  \bE \left[ 
B_1(X_{t_i}) B_0(X_{t_i}) I_{(1), i} \dt
B_1(X_{t_j}) B_2(X_{t_j}) I_{(1), j} I_{(1,1), j}
 \right] = \\
& \frac{1}{M^2 \dt}    \sum\limits_{j \in M_k}  \bE \left[
B_1^2(X_{t_j}) B_0(X_{t_j}) B_2(X_{t_j}) I_{(1), j}^2 I_{(1,1), j}
 \right]  = \\
& \frac{1}{M^2 \dt}  \sum\limits_{j \in M_k}  \bE_k \left[
B_1^2(X_{t_j}) B_0(X_{t_j}) B_2(X_{t_j}) \right] \bE \left[ I_{(1), j}^2 I_{(1,1), j}
 \right] 
  \le \frac{C\dt}{M}, 
\end{align*}
where we use Lemma 5.7.2 and Lemma 5.7.5 in \cite{Peter} to obtain the order of $\dt$. 
Other terms result in higher-order terms. Therefore, Type 2a terms are equivalent to $O(\dt/M)$.

\noindent
{\bf{Type 2b}:} \\
Consider type 2 terms with the following 3 restrictions - \\ 
(i) at least one of the integrals in each pair $I_{\alpha_q, i} I_{\alpha_l, i}$ and $I_{\alpha_r, j} I_{\alpha_m, j}$ is stochastic,\\
(ii) $q \ne l$ and $r \ne m$,\\
(iii) $n(\alpha_q, \alpha_l)$ or  $n(\alpha_r, \alpha_m)$ is even.

Clearly, condition (iii) here is complimentary 
to the condition (iii) for type 2a terms. For type 2b terms the summation over $i,j \in M_k$ 
cannot be reduced to the summation $j \in M_k$. 
Thus, we provide different types of estimates compared with type 2a terms.
In particular, we consider  
\begin{eqnarray}
&& \frac{1}{(M\dt)^2}  \labs \bE \left[\sum\limits_{i,j \in M_k} 
B_q(X_{t_i}) B_l(X_{t_i}) I_{\alpha_q, i} I_{\alpha_l, i}
B_r(X_{t_j}) B_m(X_{t_j}) I_{\alpha_r, j} I_{\alpha_m, j} 
 \right] \rabs \nonumber \le 
 \nonumber \\
 &&  \frac{1}{\dt^2}  \| B_q(X_{t_i}) B_l(X_{t_i}) B_r(X_{t_j}) B_m(X_{t_j}) \|_2 \,
 \| I_{\alpha_q, i} I_{\alpha_l, i} I_{\alpha_r, j} I_{\alpha_m, j} \|_2  = 
 \nonumber \\
 && \frac{C}{\dt^2} \| I_{\alpha_q, i} I_{\alpha_l, i} I_{\alpha_r, j} I_{\alpha_m, j} \|_2.
 \label{reduction1}
\end{eqnarray}
We would like to point out that when $i \ne j$ the norm above reduces to 
\[
\| I_{\alpha_q, i} I_{\alpha_l, i} I_{\alpha_r, j} I_{\alpha_m, j} \|_2 = 
\| I_{\alpha_q, i} I_{\alpha_l, i} \|_2 \,  \| I_{\alpha_r, j} I_{\alpha_m, j} \|_2 \quad 
\text{for~} i \ne j.
\]
The constant $C$ is finite because this generic constant corresponds to the 
norm with respect to the joint conditional distribution of $X_{t_i}$ and $X_{t_j}$, i.e.,
\begin{eqnarray*}
&& \| B_q(X_{t_i}) B_l(X_{t_i}) B_r(X_{t_j}) B_m(X_{t_j}) \|_2  = \\
&& \left( \bE \left[ \left(B_q(X_{t_i}) B_l(X_{t_i}) B_r(X_{t_j}) B_m(X_{t_j})\right)^2 | X_{t_i}, X_{t_j} \in Bin_k \right] \right)^{1/2}.
\end{eqnarray*}

Therefore, we have to compute the lowest-order terms of the form
\[
\| I_{\alpha_q, i} I_{\alpha_l, i} I_{\alpha_r, j} I_{\alpha_m, j} \|_2 ,
\]
where we use Lemma 5.7.5 in \cite{Peter} for $i = j$ and Lemma 5.7.2 in \cite{Peter} for $i \ne j$.
We would like to note that without the restriction (ii) the lowest order terms would be 
$\| I_{(1),i}^2 I_{(1),j}^2\|_2 = \dt^2$. However, with restriction (ii)  
neither $I_{(1),i}^2$ nor $I_{(1),j}^2$ are allowed in the summation.
We demonstrate here that type 2b terms are equivalent to $O(\dt^{3/2})$.

Without the loss of generality we consider the case when $n(\alpha_r, \alpha_m)$ is even and
list the lowest order terms:

\noindent
(a) When $(\alpha_q, \alpha_l) = ((1), (0))$ and $(\alpha_r, \alpha_m) = ((0), (1, 1))$
\eqref{reduction1} becomes
\[
\dt^{-2} \| I_{\alpha_q, i} I_{\alpha_l, i} I_{\alpha_r, j} I_{\alpha_m, j} \|_2 =
\dt^{-2} \dt^2 \| I_{(1), i} I_{(1,1), j} \|_2 \le  C \dt^{\frac{3}{2}}
\]
for both, $i=j$ and $i \ne j$. 
In fact, the norm above can be computed exactly in both cases since both 
$I_{(1), i}$ and $I_{(1,1), i}$ can be represented explicitly through the increment
of the Brownian motion $\Delta W_{j+1}$.

\noindent
(b) When $(\alpha_q, \alpha_l) = ((1), (0))$ and $(\alpha_r, \alpha_m) = ((1), (0, 1))$ \eqref{reduction1} becomes
\[
\dt^{-2} \| I_{\alpha_q, i} I_{\alpha_l, i} I_{\alpha_r, j} I_{\alpha_m, j} \|_2 =
\dt^{-2} \dt \| I_{(1), i} I_{(1,1), j} I_{(0,1), j} \|_2 \le  C \dt^{\frac{3}{2}}
\]
for both $i = j$ and $i \ne j$. We use the Minkowski inequality and omit higher order terms for $i = j$. Considering $i \ne j$, the stochastic integrals with different subscripts are independent, we can separate the $L^2$ norm in \eqref{reduction1} into a product of two $L^2$ norms.

\noindent
(c) When $(\alpha_q, \alpha_l) = ((1), (0))$ and $(\alpha_r, \alpha_m) = ((1), (1, 0))$ \eqref{reduction1} becomes
\[
\dt^{-2} \| I_{\alpha_q, i} I_{\alpha_l, i} I_{\alpha_r, j} I_{\alpha_m, j} \|_2 =
\dt^{-2} \dt \| I_{(1), i} I_{(1,1), j} I_{(1,0), j} \|_2 \le  C \dt^{\frac{3}{2}}
\]
for both $i = j$ and $i \ne j$. Calculation is similar to the case (b) above.

\noindent
(d) When $(\alpha_q, \alpha_l) = ((1), (1, 1))$ and $(\alpha_r, \alpha_m) = ((0), (1, 1))$, we have
\[
\dt^{-2} \| I_{\alpha_q, i} I_{\alpha_l, i} I_{\alpha_r, j} I_{\alpha_m, j} \|_2 =
\dt^{-1} \| I_{(1), i} I_{(1,1), i} I_{(1,1), j} \|_2 \le  C \dt^{\frac{3}{2}}
\]
for both $i = j$ and $i \ne j$ since all integrals can be explicitly expressed through $\dt$ and
$\Delta W_{j+1}$.

\noindent
(e) When $(\alpha_q, \alpha_l) = ((1), (1, 1))$ and $(\alpha_r, \alpha_m) = ((1), (0, 1))$, we have
\[
\dt^{-2} \| I_{\alpha_q, i} I_{\alpha_l, i} I_{\alpha_r, j} I_{\alpha_m, j} \|_2 =
\dt^{-2} \| I_{(1), i} I_{(1,1), i} I_{(1), j} I_{(0,1), j} \|_2 \le  C \dt^{\frac{3}{2}}
\]
for both $i = j$ and $i \ne j$. When $i = j$ we use  Lemma 5.7.2 in \cite{Peter} and
if  $i \ne j$ we can use properties of stochastic integrals in section \ref{sec:sec3}.

\noindent
(f) When $(\alpha_q, \alpha_l) = ((1), (1, 1))$ and $(\alpha_r, \alpha_m) = ((1), (1, 0))$, we have
\[
\dt^{-2} \| I_{\alpha_q, i} I_{\alpha_l, i} I_{\alpha_r, j} I_{\alpha_m, j} \|_2 =
\dt^{-2} \| I_{(1), i} I_{(1,1), i} I_{(1), j} I_{(1,0), j} \|_2 \le  C \dt^{\frac{3}{2}}
\]
for both $i = j$ and $i \ne j$. Calculations here are similar to the case (e) above.

\noindent
(g) When $(\alpha_q, \alpha_l) = ((1), (1, 1))$ and $(\alpha_r, \alpha_m) = ((1), (1, 1, 1))$, we have
\[
\dt^{-2} \| I_{\alpha_q, i} I_{\alpha_l, i} I_{\alpha_r, j} I_{\alpha_m, j} \|_2 =
\dt^{-2} \| I_{(1), i} I_{(1,1), i} I_{(1), j} I_{(1,1,1), j} \|_2 \le  C \dt^{\frac{3}{2}}
\]
for both $i = j$ and $i \ne j$. From 5.2.21 of \cite{Peter}, we have
$I_{(1,1,1), j} = \frac{1}{3^^21} (I^3_{(1), j} - 3\dt I_{(1), j})$. 
Therefore, all stochastic integrals can be expressed explicitly through 
$\Delta W_{j+1}$.

\noindent
(h) When $(\alpha_q, \alpha_l) = ((1), (0))$ and $(\alpha_r, \alpha_m) = ((1), (1, 1, 1))$, we have
\[
\dt^{-2} \| I_{\alpha_q, i} I_{\alpha_l, i} I_{\alpha_r, j} I_{\alpha_m, j} \|_2 =
\dt^{-1} \| I_{(1), i} I_{(1), j} I_{(1,1,1), j} \|_2 \le  C \dt^{\frac{3}{2}}
\]
for both $i = j$ and $i \ne j$. Similar to the previous case, we can express
$I_{(1,1,1), j} = \frac{1}{3^^21} (I^3_{(1), j} - 3\dt I_{(1), j})$. 
Therefore, all stochastic integrals can be expressed explicitly through 
$\Delta W_{j+1}$.

\noindent
{\bf{Type 2c}:} \\
Here we consider the case when both integrals in one pair $I_{\alpha_q, i} I_{\alpha_l, i}$ or $I_{\alpha_r, j} I_{\alpha_m, j}$ are deterministic.
Without the loss of generality we consider both integrals $I_{\alpha_q, i} I_{\alpha_l, i}$ to be deterministic.
There are only two deterministic integrals considered in the truncated Ito-Taylor expansion \eqref{ITE}, namely
$I_{(0),i} = \dt$ and $I_{(0,0),i} = \dt^2/2$.
Clearly, the lowest-order terms arise from $(\alpha_q, \alpha_l) = ((0), (0))$ (i.e., $I_{\alpha_q, i} I_{\alpha_l, i} = I_{(0),i}^2 = \dt^2$). 

We use the same approach as for type 2b terms in \eqref{reduction1}. In particular, we obtain estimate
\[
\frac{1}{(M\dt)^2}  \labs \bE \left[\sum\limits_{i,j \in M_k} 
B_q(X_{t_i}) B_l(X_{t_i}) I_{\alpha_q, i} I_{\alpha_l, i}
B_r(X_{t_j}) B_m(X_{t_j}) I_{\alpha_r, j} I_{\alpha_m, j} 
 \right] \rabs \le
C \| I_{\alpha_r, j} I_{\alpha_m, j} \|_2.
\]
Here we list lowest-order terms.

\noindent
(a) When $(\alpha_q, \alpha_l)= ((0), (0))$ and $(\alpha_r, \alpha_m)=((1), (0))$ 
(or $((0), (1))$) we obtain
\[
\| I_{\alpha_r, j} I_{\alpha_m, j} \|_2 =
 \dt \| I_{(1), j} \|_2 = 
 \dt^{3/2}. 
\]

\noindent
(b) When $(\alpha_q, \alpha_l) = ((0), (0))$ and $(\alpha_r, \alpha_m)=((1), (1, 1))$ (or $((1, 1), (1))$) we obtain
\[
\| I_{\alpha_r, j} I_{\alpha_m, j} \|_2 = 
\| I_{(1), j} I_{(1,1), j} \|_2 \le
C \dt^{3/2} 
\]
where we used \eqref{I11}.
All other combinations of stochastic integrals yield terms of higher order.
Therefore, Type 2c terms are equivalent to $O(\dt^{3/2})$.

\noindent
{\bf{Type 2d}:} \\
Consider terms with $q = l$ or $r = m$. Without loss of generality we consider the case $q=l$. We would like to remind that type 2 terms are
computed under the restriction $q \times l \ne 1$, which means that $\alpha_q = \alpha_l = (1)$ does not
occur for type 2d terms.
Stochastic integrals which yield the lowest-order terms are $q=l=2$ or $I_{(1,1),i}^2 \sim \bE \left[ \Delta W_{i+1}^4\right] \sim \dt^2$.

Here we use the same approach as for type 2b terms in \eqref{reduction1}. In particular, we write
\[
\frac{1}{(M\dt)^2}  \labs \bE \left[\sum\limits_{i,j \in M_k} 
B_q(X_{t_i}) B_l(X_{t_i}) I_{\alpha_q, i} I_{\alpha_l, i}
B_r(X_{t_j}) B_m(X_{t_j}) I_{\alpha_r, j} I_{\alpha_m, j} 
 \right] \rabs  
\le  \frac{C}{\dt^2} \| I_{(1,1),i}^2 I_{\alpha_r, j} I_{\alpha_m, j} \|_2.
\]
The lowest-order terms arise from combination of indexes $(\alpha_r, \alpha_m) = ((1), (0))$ and $(\alpha_r, \alpha_m) = ((1), (1,1))$ and we treat these two cases next.

\noindent
(a)  When $(\alpha_r, \alpha_m) = ((1), (0))$ we obtain
\[
C\dt^{-2} \| I_{(1,1),i}^2 I_{\alpha_r, j} I_{\alpha_m, j} \|_2 = 
C\dt^{-2} \dt \| I_{(1,1),i}^2 I_{(1), j} \|_2 \le
C \dt^{3/2}
\]
for both $i=j$ and $i \ne j$ where we used \eqref{I11} and the Minkowski inequality.

\noindent
(b)  When $(\alpha_r, \alpha_m) = ((1), (1,1))$ we obtain
\[
C\dt^{-2} \| I_{(1,1),i}^2 I_{\alpha_r, j} I_{\alpha_m, j} \|_2 = 
C\dt^{-2} \dt \| I_{(1,1),i}^2 I_{(1), j} I_{(1,1),j} \|_2 \le
C \dt^{3/2}
\]
for both $i = j$ and $i \ne j$ where we used \eqref{I11} and the Minkowski inequality.
All other combinations of integrals yield terms of higher order.
Therefore, Type 2d terms are equivalent to $O(\dt^{3/2})$.

\noindent
{\bf{Type 3}:} 
Finally, we consider Type 3 terms and using \eqref{I11} we obtain
\begin{eqnarray*}
&& \frac{2}{(M\dt)^2} \bE \left[  \left( \sum\limits_{i,j \in M_k}
\left( D^2(X_{t_i}) I_{(1), i}^2 - D^2(x_k) \dt \right)
\sum\limits_{\substack{q, l=0 \\q \times l \ne 1}}^6 B_q(X_{t_j}) B_l(X_{t_j}) I_{\alpha_q, j} I_{\alpha_l, j} \right) \right]  = \\
&& \frac{2}{(M\dt)^2} \bE \left[  \left(\sum\limits_{i,j \in M_k}
 D^2(X_{t_i})  \left(I_{(1),i}^2 - \dt \right)
\sum\limits_{\substack{q, l=0 \\ q \times l \ne 1}}^6 B_q(X_{t_j}) B_l(X_{t_j}) I_{\alpha_q, j} I_{\alpha_l, j} \right)
\right] + \\
&&  \frac{2}{(M\dt)^2} \bE \left[  \left(\sum\limits_{i,j \in M_k}
\dt  \left( D^2(X_{t_i}) -  D^2(x_k) \right)
\sum\limits_{\substack{q, l=0 \\ q \times l \ne 1}}^6 B_q(X_{t_j}) B_l(X_{t_j}) I_{\alpha_q, j} I_{\alpha_l, j} \right)\right]  = \\
&&  \underbrace{\frac{4}{(M\dt)^2} \bE \left[  \left(\sum\limits_{i,j \in M_k} 
 D^2(X_{t_i}) I_{(1, 1),i}
\sum\limits_{\substack{q, l=0 \\ q \times l \ne 1}}^6 B_q(X_{t_j}) B_l(X_{t_j}) I_{\alpha_q, j} I_{\alpha_l, j} \right)
\right]}_{Type~3a}  + \\
&&  \underbrace{\frac{2}{M^2 \dt} \bE \left[  \left(\sum\limits_{i,j \in M_k} 
 \left( D^2(X_{t_i}) -  D^2(x_k) \right)
\sum\limits_{\substack{q, l=0 \\ q \times l \ne 1}}^6 B_q(X_{t_j}) B_l(X_{t_j}) I_{\alpha_q, j} I_{\alpha_l, j} \right)\right]}_{Type~3b},
\end{eqnarray*}
where we used \eqref{I11}.

\noindent
{\bf{Type 3a}:} 
For the first term in type 3a, the lowest-order terms arise from $q = 1$, $l = 0$ and $q = 1$, $l = 2$ which corresponds to 
$(\alpha_q, \alpha_l) = ((1), (0))$ and $(\alpha_q, \alpha_l) = ((1), (1,1))$, respectively.

\noindent
(a) Consider $q = 1$, $l = 0$ first. 
Then using the same argument as for type 2a we can show that
\begin{align*}
& \frac{4}{(M\dt)^2} \bE \left[  \sum\limits_{i,j \in M_k}
 D^2(X_{t_i})  I_{(1,1), i}
 B_1(X_{t_j}) B_0(X_{t_j}) I_{(1), j} \dt 
 \right]   = \\
& \frac{4}{M^2 \dt} \bE \left[  \sum\limits_{i \in M_k}
 D^2(X_{t_i})  I_{(1,1), i}
 B_1(X_{t_i}) B_0(X_{t_i}) I_{(1), i} 
 \right]    = \\
& \frac{4}{M^2 \dt} \sum\limits_{i \in M_k}
 \bE_k \left[ D^2(X_{t_i}) 
 B_1(X_{t_i}) B_0(X_{t_i})  \right] \bE \left[ I_{(1,1), i} I_{(1), i} 
 \right]   = 0.
\end{align*}

\noindent
(b) Next, consider $q = 1$, $l = 2$. Then using the same argument as above we obtain
\begin{align*}
& \frac{4}{(M\dt)^2} \bE \left[  \sum\limits_{i,j \in M_k}
 D^2(X_{t_i})  I_{(1,1), i}
 B_1(X_{t_j}) B_2(X_{t_j}) I_{(1), j} I_{(1,1), j}
 \right]    = \\
& \frac{4}{(M\dt)^2} \bE \left[  \sum\limits_{i \in M_k}
 D^2(X_{t_i})  I_{(1,1), i}
 B_1(X_{t_i}) B_2(X_{t_i}) I_{(1), i} I_{(1,1), i}
 \right]   = \\
 & \frac{4}{(M\dt)^2} \sum\limits_{i \in M_k}
\bE_k \left[ D^2(X_{t_i}) 
 A_1(X_{t_i}) A_5(X_{t_i}) \right] \bE \left[ I_{(1), i} I_{(1,1), i}^2
 \right] = 0.
\end{align*}

\noindent
(c) The next order terms appear due to combinations of indexes which correspond to
$(\alpha_q, \alpha_l) = ((1), (0,1))$, $(\alpha_q, \alpha_l) = ((1), (1,0))$, and $(\alpha_q, \alpha_l) = ((1), (1,1,1))$.
In these cases we cannot apply argument used perviously in (a) and (b) since
 $\bE\left[ I_{(1), i} I_{(1,0), i} \right] \ne 0$, $\bE\left[ I_{(1), i} I_{(0,1), i} \right] \ne 0$, and 
$\bE\left[ I_{(1), i} I_{(1,1,1), i} \right] \ne 0$.
Therefore, we first consider the case $(\alpha_q, \alpha_l) = ((1), (0,1))$ and
proceed as in \eqref{reduction1} to obtain 
\begin{align*}
& \frac{4}{(M\dt)^2} \bE \left[  \sum\limits_{i,j \in M_k}
 D^2(X_{t_i})  I_{(1,1), i}
 B_1(X_{t_j}) B_4(X_{t_j}) I_{(1), j} I_{(1,0), j}
 \right] \le    \\
 & \frac{4}{\dt^2} \| D^2(X_{t_i})  B_1(X_{t_j}) B_4(X_{t_j}) \|_2 \,
 \| I_{(1,1), i} I_{(1), j} I_{(1,0), j} \|_2 \le 
 C \dt.
\end{align*}
A similar argument can be applied to $(\alpha_q, \alpha_l) = ((1), (0,1))$ and $(\alpha_q, \alpha_l) = ((1), (1,1,1))$ to yield the same bound $O(\dt)$.

\noindent
(d) We would like to point out that the combination of indexes
$(\alpha_q, \alpha_l) = ((0), (1,1))$ yields a higher-order term because in this case we can use the argument similar (a) and (b) to obtain
\begin{align*}
& \frac{4}{(M\dt)^2} \bE \left[  \sum\limits_{i,j \in M_k}
 D^2(X_{t_i})  I_{(1,1), i}
 B_0(X_{t_j}) B_2(X_{t_j}) \dt I_{(1,1), j}
 \right]    = \\
 & \frac{4}{M^2 \dt } \bE \left[  \sum\limits_{i,j \in M_k}
 D^2(X_{t_i})  I_{(1,1), i}
 B_0(X_{t_j}) B_2(X_{t_j}) I_{(1,1), j}
 \right]  = \\
 & \frac{4}{M^2 \dt} \bE \left[  \sum\limits_{i \in M_k}
 D^2(X_{t_i}) 
 B_0(X_{t_i}) B_2(X_{t_i}) I_{(1,1), i}^2
 \right]  \le \\
 & \frac{4}{M^2 \dt}  \sum\limits_{i \in M_k}
 \bE_k \left[ D^2(X_{t_i}) B_0(X_{t_i}) B_2(X_{t_i})  \right] \,
 \bE \left[   I_{(1,1), i}^2 \right]
 \le \frac{C \dt}{M}.
\end{align*}

\noindent
{\bf{Type 3b}:} 
\begin{align*}
& \frac{2}{(M\dt)^2} \left\lvert \bE \left[  \sum\limits_{i,j \in M_k}  
\dt  \left( D^2(X_{t_i}) -  D^2(x_k) \right)
\sum\limits_{\substack{q, l=0 \\ q \times l \ne 1}}^6 B_q(X_{t_j})B_l(X_{t_j}) I_{\alpha_q, j} I_{\alpha_l, j} 
 \right] \right\rvert \le \\
&  \frac{2K_d \dx}{\dt} 
\sum\limits_{\substack{q, l=0 \\ q \times l \ne 1}}^6
\bE \left[ \left\lvert 
B_q(X_{t_j}) B_l(X_{t_j}) I_{\alpha_q, j} I_{\alpha_l, j} \right \rvert  \right]  = \\
& \frac{2K_B \dx}{\dt} 
\sum\limits_{\substack{q, l=0 \\ q \times l \ne 1}}^6
\bE_k \left[ \left\lvert 
B_q(x) B_l(x)  \right\rvert \right]
\| I_{\alpha_q, j} I_{\alpha_l, j} \|_1.
\end{align*}

Lowest-order terms arise from 
$q = 1$, $l = 0$ or $q = 1$, $l = 5$. Consider $q = 1$, $l = 0$. 
Then $\| I_{(1), j} I_{(0), j} \|_1 = \dt \| I_{(1), j} \|_1 = O(\dt^{3/2})$ and
\begin{equation*}
 \frac{2K_d \dx}{\dt} 
\bE_k \left[ \left\lvert 
B_1(x) B_0(x)  \right\rvert \right]
\| I_{(1), j} I_{(0), j} \|_1 \le C \dx \sqrt{\dt} .
\end{equation*}
One can also show that $\| I_{(1), j} I_{(1,1), j} \|_1 \sim \bE [|\Delta W_{j+1}^3|] = O(\dt^{3/2})$ which 
yields a similar bound for $q = 1$, $l = 2$.

\section*{Acknowledgements}
This research has been partially supported by grants NSF DMS-1620278 and ONR N00014-17-1-2845.


\end{document}